\documentclass[11pt]{article}
\usepackage{amsfonts}
\usepackage{mathrsfs}
\usepackage{amsmath,amsfonts,mathrsfs,amssymb,color}
\usepackage{indentfirst}
\usepackage{graphicx}
\usepackage{float}
\usepackage{subfigure}

\numberwithin{equation}{section}

\newtheorem{theo.}{\quad\, Theorem}[section]
\newtheorem{defi.}{\quad\, Definition}[section]
\newtheorem{lemm.}{\quad\, Lemma}[section]
\newtheorem{coro.}{\quad\, Corollary}[section]

\begin{document}

\title{On the well-posedness and efficient approximation for the classical Melan equation in suspension bridges$^*$ }
\author{Jinxiang Wang$^{\ast}$%, \ \ \ \ Da-Bin Wang
\\
 {\small  Department of Applied Mathematics, Lanzhou University of Technology, Lanzhou, P.R. Chin}\\
}
\date{} \maketitle
\footnote[0]{E-mail address: wjx19860420@163.com(Jinxiang Wang),  %wangdb96@163.com (Da-Bin Wang) \
%Telephone: 86-931-7971297
% $^{\ddagger} $xuling$_{-}$216@yahoo.cn
} \footnote[0] {$^*$Corresponding author: Jinxiang Wang.\   }

\begin{abstract}
\baselineskip 18pt
The classical Melan equation modeling suspension bridges is considered. We first study the explicit expression and global properties of the analytical solution for the simplified ``less stiff'' model, based on which we derive an a priori estimate for the original classical Melan equation and establish, under an explicit load condition, that it has a unique solution with nonnegative integral under downward live loads, thereby showing the uniqueness of the corresponding deflection curve in the engineering setting. We also develop an efficient iterative approximation method by taking the solution of the simplified ``less stiff'' model as the first iterate, and prove its geometric convergence with explicit error estimates. The applicability and computational efficiency of the method are demonstrated through calculations for two actual bridges, which also quantify the influence of the nonlinear nonlocal term on the solution and clarify the relationship between the simplified and original models. Several engineering observations are verified and explained, and some related open problems are suggested.
 \end{abstract}


\vskip 3mm


\baselineskip 22pt

\section{Introduction}

%Mathematical models
The deflection theory for suspension bridges was
originally developed by the Austrian engineer Josef Melan in monograph [1] whose first edition goes back to 1888. As the culmination of a series of suspension bridge theories that began with Navier's unstiffened bridge theory [2] in 1823 and included the Rankine's elastic theory of the stiffened suspension bridge(see [3]), Melan's deflection theory is acknowledged as a milestone contribution to the understanding of suspension bridges. A historical overview of the methods leading to deflection theory can be found in the fascinating survey of Buonopane and Billington [4] and Pugsley [5]. Since the famous Manhattan Bridge was first designed and built using deflection theory in 1909, the accuracy of deflection theory in describing the structural characteristics of suspension bridges has been continuously verified and recognized. Combined with the great economic benefits brought by deflection theory in the design of the deck stiffening, more and more long-span suspension bridges have been designed and built under the guidance of the theory.  % ( Since the famous Manhattan Bridge, which was first designed and constructed using deflection theory in 1909, the accuracy of deflection theory in describing the structural features of suspension Bridges has been continuously verified. This, combined with the economic benefits brought by this theory, has led to the design and construction of an increasing number of long-span suspension Bridges under its guidance.)
%In 1929, Steinman included the deflection theory in the second edition of his classic Practical Treatise [5]. Since then, more and more...
Today, although the engineering design of suspension bridges is often performed using computer methods that employ discrete structural elements and are grounded in finite displacement theory (see, e.g., [6-7]), %although advances in computational methods have prompted several new approaches to analysis of suspension-bridges, Modern suspension bridges are typically analyzed using computer programs with nonlinear analysis capabilities based on finite-element formulations.
the deflection theory is still used as the main analytical method in the preliminary design and practical calculation of suspension bridges (see, e.g., [8-19]), which enables designers to understand the influence of key parameters and the basic behavior of structure quickly, and to validate the complex simulations.

%From the very beginning, the core concern of suspension bridge theory is the vertical deflection of the bridge deck under an applied live load, and the corresponding mathematical models usually are one-dimensional in nature. %The core question of the theory is primarily concerned with vertical deflection of the bridge deck under an applied live load, and the corresponding Mathematical models usually are one-dimensional in nature.
In deflection theory,  Melan [1] models a suspension bridge as an elastic beam (the deck) suspended to a sustaining cable (see Figure 1 below).
\begin{figure}[H]
\centering
\includegraphics[width=1.0\linewidth]{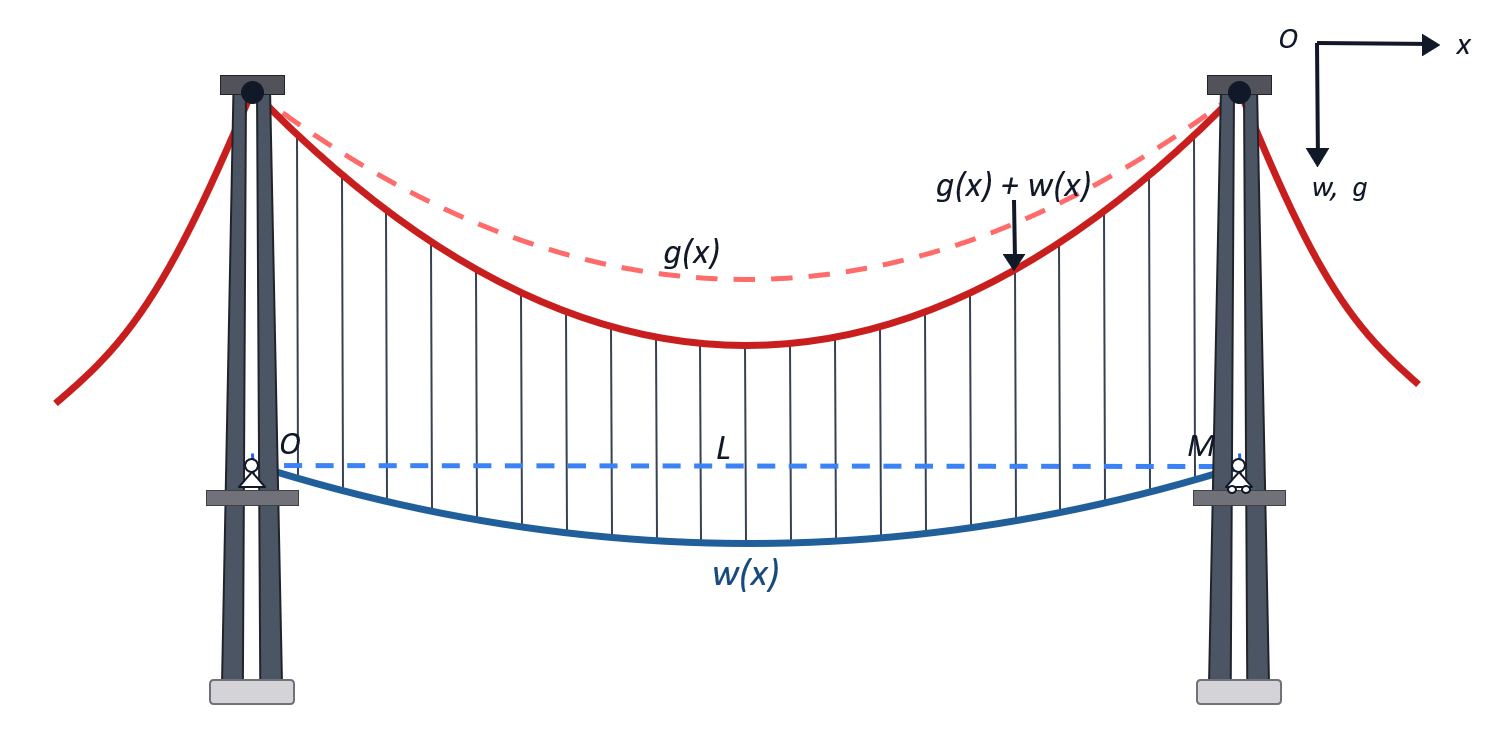}
\caption{ Beam (blue) sustained by a cable (red) through
parallel hangers.}
\end{figure}
\noindent
To interpret the core problems such as the vertical deflection of the bridge deck under live load and the moments and shears in a suspension bridge, Melan [1, p.77] proposed a fourth-order equation to describe the deck deflection, which is expressed as:
$$
EI w''''(x)-(H+h(w))w''(x)+\frac{q}{H}h(w)=P(x),\ \ \ \ x\in (0,L),
\eqno (1.1)
$$
where $w=w(x)$ denotes the vertical displacement of the beam representing the deck (positive downward), $L$ is the distance between the two towers and it is also the span of the beam hinged at the two endpoints. This means that the boundary conditions for equation (1.1) are $$
w(0)=w(L)=w''(0)=w''(L)=0.\eqno (1.2)
$$
$E$ and $I$ in (1.1) are respectively the elastic modulus of the material composing the deck and the moment of inertia of the cross section so that $EI$ is the flexural rigidity, $H$ is the horizontal tension of the cable when subject to the dead load $q$, and $h(w)$ is a nonlocal term under the integral sign (detailed in below) representing the additional tension in the cable due to the live load $P=P(x)$. %.   The dead load $q$ includes the weights of the cable, of the hangers, and of the deck.

Melan's equation (1.1) is called {\it the fundamental equation of the theory of suspension bridge} by Von K\'arm\'an and Biot in the monograph [20], in which a specific approximate expression of the nonlocal term $h(w)$ is also given (see (5.14) in [20]) as follows:
$$
h(w)=\frac{E_{c}A_{c}}{L_{c}}\frac{q}{H}\int_0^L w(x)dx,\eqno (1.3)
$$
where $A_{c}$ denotes the cross-sectional area of the cable and $E_{c}$ denotes the modulus of elasticity of the material, $L_{c}$ is the length of the cable subject to the dead load $q$. The above approximate formula and the corresponding specific form of Melan equation (substituting (1.3) into (1.1)):
$$
EI w''''(x)-(H+\frac{E_{c}A_{c}}{L_{c}}\frac{q}{H}\int_0^L w(x)dx)w''(x)+\frac{q^{2}}{H^{2}}\frac{E_{c}A_{c}}{L_{c}}\int_0^L w(x)dx=P(x)\eqno (1.4)
$$
given by [20] % (see Appendix in Section 8 for detailed derivation of the equation)  %have been recognized and
have since been widely adopted as the classical governing equation for deflection theory in engineering literature and textbooks due to their conciseness.%(see (11.17) in Timoshenko and Young [22];).%see [6-16].      are still widely used in textbooks and related literatures, see...      In the monograph [17] by Von K′arm′an and Biot, (1.1) is called {\it the fundamental equation of the theory of suspension bridge}, and this term has been recognized and used by subsequent literature and textbooks until today. %For the nonlocal item $h(w)$ in (1.1), the most common expression for nonlocal terms are as followsthere are three different expressions at present, (see [17, 18, 20]). In (1.1), the most common expression for the nonlocal term $h(w)$ is as follow: ...$h(w)$ is a nonlocal term under the integral sign representing the additional tension in the cable due to the live load $p=p(x).$ There are three different expressions for the nonlocal term $h(w)$ in the existing literature which will be detailed below in Section2.

%The deflection theory of Melan and extensions of it are the only long-established methods of performing static analyses of suspension bridges with any degree of accuracy.
From a mathematical point of view, the term $h(w)$ depending on the deflections $w$ with determination (1.3) %under the integral sign
 makes (1.1) a nonlinear nonlocal equation, which brings great acknowledged difficulties to the qualitative research of analytical solutions. For this reason, since the deflection theory of suspension bridges was
established more than 100 years ago, several simplifications and approximate methods have been developed in engineering literature for analysing and solving the fundamental equations (1.1) and (1.4). One of the most commonly used approximate solution methods for (1.1) and (1.4) is the trial-and-error method (see Steinmann [21]): a value for $h$ is first assumed to solve (1.1) and (1.2); the obtained expression for $w$ is then substituted into the expression for $h(w)$ similar to (1.3) for verification. If the outcome does not match the assumed value, $h$ is updated and the calculation is repeated until a solution satisfying the required accuracy is obtained. This method was once referred to as the ``exact method,'' but the procedure was tedious before the advent of computers and lacks rigorous theoretical analysis. %In [21] and [22], the method of trigonometric series was proposed: first treats $h$ as a constant, employs trigonometric series to derive a solution for $w$ expressed as a sine series containing $h$, and then determines $h$ by combining the trial-and-error method with interpolation technique.
For other simplifications and approximate methods to fundamental equations (1.1) and (1.4) in engineering literature, refer to, for example, the method of influence lines in [22] and [23], the linearized deflection theory in [24], the method of trigonometric series in [20] and [25], the substitutional beam method in [10], the gravity stiffness method in [26], and the stiffness matrix method based on Laplace transformation in [27]. As can be seen from [20-27] discussed above and the later related works [13,14], the simplifications and approximate methods in existing engineering literature generally lack rigorous analysis and proof. In particular, the relationship between approximate and true solutions is rarely addressed, and convergence proofs are entirely absent.
In [28-29], Semper considered the ``less stiff''
case of the fundamental equation (1.1) and (1.4): If the live load is relatively small compared to the dead load, $H$ will be much greater than $h(w)$, then $H+h(w)\approx H$ %. For such case, Semper assumed that H + h(w) ≈ H,
 and the equation (1.4) can be reduced to the following:

$$
EI w''''(x)-Hw''(x)+\frac{q^{2}}{H^{2}}\frac{E_{c}A_{c}}{L_{c}}\int_0^L w(x)dx=P(x),\ \ \ \ x\in (0,L).
\eqno (1.5)
$$
Semper developed a finite element method for (1.5),(1.2) in [28], and the convergence of a Galerkin approximation for this linear integro-differential problem was further analyzed in [29]. Numerical comparisons in [28] suggested that (1.5) may serve as an approximation to the nonlinear nonlocal equation (1.4) in the examples considered; however, no rigorous relationship or quantitative error estimate between the two models was established. %The linearized deflection theory given by (1.5) is recognized as a satisfactory approximation to (1.4) with specific small live loads,%In a subsequent work [21], Semper discussed the convergence of the iteration algorithm to the approximate finite-element solution of (1.4),(1.2).
In [30], Gauss-Jacobi iterative technique based on a spectral method is given to get the approximate solution of (1.5), (1.2), and the iterative results similar to [28-29] are obtained. %that is, the non-linear term in equation (1.4) does not seem to make appreciable contribution.
For other studies on numerical approximation methods for (1.5), (1.2), refer to the Legendre-Galerkin spectral approximation in [31], the modified Arithmetic Mean iterative method in [32], and the Optimal Homotopic Asymptotic Method in [33].

Although numerical experiments in [28] suggest that the nonlinear nonlocal term in equation (1.4) seems to have only a limited influence in the examples considered, as the nonlinear structural characteristics of the flexible structure of suspension bridge have been discovered and discussed (see, e.g., [34-37]), the qualitative analysis on the original nonlinear nonlocal fundamental equation (1.4) has become an unavoidable problem.
%As the nonlinear structural characteristics of the flexible structure of suspension bridge are constantly discovered and discussed(see, e.g., [25, 26, 27, 28]), the qualitative research on the original nonlinear basic equation which has not been simplified has become an inevitable problem. However, the nonlinear structural behavior of suspension bridges is by now well established; see, e.g., [3, 10, 12, 15, 21].% and (1.5),(1.2) has been used to analyze the structure of Bridges in real life, the problem of ()
%It is worth noticing that, in the above mentioned engineering literature [6-19] and numerical analysis works [20-24], qualitative problems such as the analytical solutions of simplified ``less stiff'' model (1.5),(1.2) and the existence of solutions for the fundamental equation (1.1),(1.2) and its classical form (1.4),(1.2)  have not been discussed.
In 2014, %inspired by some newly discovered nonlinear structural characteristics of suspension bridges(see, e.g., [25, 26, 27, 28]),
Gazzola et al. [38] paid special attention to the nonlocal term $h(w)$ in the fundamental equation (1.1) and carried out a pioneering qualitative analysis of problems (1.1),(1.2) and (1.4),(1.2). Besides the common expression (1.3), they considered two other possible forms of
$h(w)$: one due to Timoshenko (see (3.9) in [38]) and another of their own (see (3.10) in [38]), and analyzed the differences among them. They also established existence results and, under stronger smallness assumptions on $L$ and $P$, an existence and uniqueness result within a class of small solutions; see [38, Theorems 5.1 and 5.4].
% ( these results are valid only for sufficiently small parameters $L$ and $p$ and, moreover, apply exclusively to small solutions )
%In addition to the common expression (1.3), two other possible forms of $h(w)$ were proposed and the differences among them were analyzed in [36]. To prove the existence and uniqueness results of solutions for (1.1), (1.2), Gazzola et al. first discussed the a priori bound for the corresponding simple constant-coefficient linear problem, % and construct a closed ball under some general(Boundedness and continuity conditions) assumptions on $h$. and consequently, under some general assumptions on $h$ the second derivative term in (1.1) was treated as a whole and then the main results were proved by using fixed point theorems when $L$ and $p$ are small; see [36, Theorem 5.1 and Theorem 5.4]. It should be emphasized that the existence and uniqueness results in [36] hold only for small solutions.
%A counterexample given by Gazzola et al. shows that: besides the small solution, there may exist additional large solutions, that is, the problem (1.1), (1.2) and (1.4), (1.2) appear to be ill-posed.
In [38], an iterative procedure based upon the solution of the simple constant-coefficient linear problem is also proposed to approximate the solutions of (1.1), (1.2) and (1.4), (1.2).  Some numerical experiments are carried out, and the results show that for some values of the parameters in (1.1) and (1.4), the iterative sequence seems to exhibit a certain quasi-monotonic behavior and admits a unique stable fixed point, whereas for some other parameters (particularly those within the practical range of bridge engineering) it appears to diverge and to be quite unstable. %At the end of the [36], Gazzola et al. proposed several open problems concerning the well-posedness of (1.1), (1.2) and (1.4), (1.2), the convergence conditions of the iterative scheme, and the possibility of developing more robust algorithms capable of handling both stable and unstable cases as well as detecting multiple fixed points.
%Noticing that, the monotonicity and convergence of the sequences observed above under certain parameters are not strictly proved,
At the end of [38], Gazzola et al. proposed several open problems: Under what conditions can the well-posedness of (1.1), (1.2) and (1.4), (1.2) be ensured? Can the solution be approximated by a suitable constructive sequence? Under which assumptions on the parameters is the iterative scheme convergent? Are there better algorithms able to manage both the stable and unstable cases? Can these algorithms detect multiple fixed points? These open problems were also listed in the subsequent monograph by Gazzola [39].
%including: identifying the conditions that ensure the well-posedness of the boundary value problems (1.1)–(1.2) and (1.4)–(1.2); investigating whether solutions can be approximated by a suitable constructive sequence; determining the parameter assumptions required for the convergence of the proposed iterative scheme; developing more robust algorithms capable of handling both stable and unstable cases; and examining whether such algorithms can detect the presence of multiple fixed points.

Recently, [40] derived the explicit expression of the unique analytical solution for the simplified ``less stiff'' model (1.5), (1.2) and, by establishing its uniform positivity within a certain range of the parameters, the author developed a monotone iterative technique based on lower and upper solutions to prove the existence, uniqueness, and approximability of the solution for the original classical fundamental equation (1.4) and (1.2). In this method, two monotone iterative sequences were constructed using the solution of the specially modified model (1.5), (1.2), which converge uniformly to the maximum and minimum solutions, respectively (or, under suitable conditions on the parameters, to the unique solution), of problem (1.4), (1.2) between the lower and upper solutions. Although the proposed method provides a rigorously analyzed and provably convergent approximation for the fundamental equation (1.4) and (1.2) and can thus be regarded as a partial answer to the open problems regarding the approximability of solutions proposed by Gazzola et al. in [38] and [39], the condition it imposes on the parameters appears overly restrictive in actual bridge validations---it is applicable only to bridges with short main span $L$ and high deck stiffness $EI$. Furthermore, the resulting solution $w$ of (1.4), (1.2) is itself subject to the constraint of positivity---that is, for admissible $P \geq 0$, the solution $w$ satisfies $w(x)\geq 0$ and $w''(x) \leq 0$, $\forall x \in (0, L)$ (see Theorems 4.1, 5.1 and Remark 4.1 in [40] for details). This implies that the deflection curve of short-span, high deck-stiffness suspension bridges is of single-curvature type without inflection points, which is consistent with the deformation pattern of a simply supported beam (see [41], p.299). This restrictive result aligns with the engineering experience that short-span, high deck-stiffness suspension bridges can be approximated by the Rankine's elastic theory (see [3]), which itself primarily resembles a simple beam solution (see [4], p.974).

%However, it is well known that the core advantage of suspension bridges lies in their long-span capability. Excessively high deck stiffness will significantly increase construction costs, thereby undermining their economic viability for large-span scenarios. Moreover, short-span suspension bridges with high deck stiffness, similar to the single-curvature bending behavior of simply supported beams, represent an idealized structural form with low material utilization efficiency. In engineering practice dominated by long-span suspension bridges, economically efficient medium-to-long span bridges typically achieve optimal internal force distribution and enhanced spanning capacity by designing and utilizing an appropriate degree of double-curvature bending.

However, it is well known that the core advantage of suspension bridges
lies in their long-span capability.  Excessively high deck stiffness may make the bridge response more beam-like, with a
predominantly single-curvature deflection shape, but it would also
significantly increase the structural weight and construction cost, thereby
reducing the economic advantage of suspension bridges, especially for large
spans.  In contrast, in practical engineering, medium- to long-span
suspension bridges, which are far more common, usually employ relatively
flexible stiffening girders. Their deflection curves may exhibit inflection points under live loads, corresponding to double-curvature bending, or more generally multi-curvature bending, which may help reduce peak bending effects and improve material utilization and spanning efficiency; see, e.g.,
[4,5,8,13,14,42].

Now, several natural questions arise. For medium- to long-span suspension
bridges commonly encountered in engineering practice, can one prove the
existence and uniqueness of solutions to the original classical Melan problem
(1.4),(1.2), particularly for physically meaningful solutions that may
exhibit inflection points (i.e., double- or multi-curvature bending)?
Moreover, can such solutions be obtained by an approximation method with
proven convergence? Furthermore, can the solution of the simplified "less stiff" model (1.5), (1.2) be directly used to approximate that of the original classical Melan problem (1.4), (1.2) without restricting the main span to be short? And can the differences and intrinsic relationship between the solutions of the simplified model and the original nonlinear nonlocal model be established theoretically?

In this paper, we address these questions by proving that, under an explicit condition on the downward live load and without imposing a separate smallness assumption on the main span length $L$, the classical Melan problem (1.4),(1.2) has a unique solution with nonnegative integral, so that the corresponding deflection curve, possibly with inflection points, is uniquely determined in the engineering setting. Under the same condition, we also develop an iterative approximation scheme by taking the solution of the simplified ``less stiff'' model (1.5),(1.2) as the first iterate and prove its geometric convergence with explicit error estimates. The actual-bridge examples in Section 4 show that this condition is satisfied for realistic small and moderate downward live loads representative of ordinary service traffic. The theoretical error analysis and numerical results show that the first iterate, namely the solution of the simplified model (1.5),(1.2), already provides a practically acceptable approximation to the solution of the original nonlinear nonlocal model (1.4),(1.2), while the second iterate yields a considerably more accurate one. In particular, the numerical examples suggest that the method may remain effective even for very heavy live loads beyond the explicit condition.

Our existence and uniqueness result extends the results of Gazzola et al.
[38] for small solutions of the classical Melan equation (1.4),(1.2) under
 small parameters \(L\) and \(P\) (see [38, Theorems 5.1 and
5.4]), as well as the result of Wang [40] for suspension bridges with short
main span length \(L\) and high deck stiffness \(EI\) (see [40, Theorem
6.1 and Remark 6.2]). In this sense, our result gives a substantial partial answer to the well-posedness problem for the classical Melan equation (1.4),(1.2) raised in [38,39] by establishing the existence and uniqueness of a solution in the physically meaningful class under realistic small and moderate downward live loads satisfying the explicit condition, without imposing separate smallness assumptions on $L$ or on the other structural parameters. Our proposed iterative method, directly based on the solution of the simplified ``less stiff'' model (1.5), (1.2), is, to the best of our knowledge, the first approximation algorithm in more than a century of suspension-bridge deflection theory with a rigorous convergence proof for practical suspension bridges of common span ranges under ordinary live loads. It also gives a substantial partial answer to the approximability problems for the classical Melan equation (1.4),(1.2) raised in [38,39]. Moreover, the theoretical error estimates and numerical results establish a quantitative connection between the solution of the simplified linear nonlocal model and that of the original nonlinear nonlocal model, thereby providing rigorous support for the earlier numerical observations reflected in [28] that the simplified ``less stiff'' model can serve as an acceptable approximation to the classical Melan model.

%Moreover, we provide a detailed error analysis between the solution of the simplified model (1.5), (1.2) and that of the original nonlinear model (1.4), (1.2), thereby offering a rigorous proof and quantitative error estimate for the numerical observation by Semper [] that the simplified model appears to approximate the original nonlinear model.

Our main tool is the Banach fixed point theorem. Early discussions of
its application to boundary value problems, including both global
contraction and contraction on a set (e.g., a ball) for second-order
problems, can be found in the monograph by Bailey et al. [43, Ch.3,
pp. 27--46] and the references therein. In 2017, Dang et al. [44] applied
the Banach fixed point theorem to nonlinear fourth-order boundary value
problems. Different from the conventional approach of defining an integral
operator on the unknown solution, they skillfully defined an operator for
the right-hand-side function from \(C[0,1]\) to \(C[0,1]\), based on the
solution of the associated basic fourth-order linear problem. Under a
Lipschitz condition on the nonlinear term with respect to the unknown
solution and its derivative within a suitably constructed bounded region,
they proved the contraction of the operator on a closed ball in \(C[0,1]\),
yielding a local existence and uniqueness result together with iterative
approximation. This approach was later extended to other beam equations;
see, e.g., [45--47]. Motivated by these works, we apply this idea to problem (1.4),(1.2) in the
present paper. Unlike [44--47], where the operator is defined through the
solution of the basic fourth-order linear boundary value problem governed
only by \(Tw=w''''\) with the corresponding boundary conditions, the
operator in the present paper is constructed from the analytical solution
and the associated integral properties obtained for the simplified
``less stiff'' model, namely the integro-differential boundary value
problem (1.5),(1.2). More importantly, we exploit the integral properties of the linear nonlocal problem to derive an a priori estimate for the original nonlinear problem (1.4),(1.2), showing that each of its solutions with nonnegative integral corresponds to a fixed point lying in a suitably constructed closed ball, and we prove that the fixed-point operator is contractive on this ball. Hence, the Banach fixed point theorem yields global uniqueness within the physically meaningful class rather than merely a local existence and uniqueness result. In this sense, the present
work extends and develops the method used in [44--47]. Furthermore, to
accommodate realistic live loads, which are naturally assumed to be
integrable and may be piecewise distributed, the present operator is
defined on the more general space \(L^{1}[0,L]\), rather than on
\(C[0,L]\) as in the previous works [44--47].

The rest of this paper is arranged as follows. In Section 2, we provide a
direct derivation of the explicit analytical solution of the simplified
``less stiff'' model (1.5),(1.2) and discuss its global properties. In
Section 3, using the solution of the simplified model and its integral
properties, we derive a priori estimates for the original nonlinear model and
establish the existence and uniqueness results for problem (1.4),(1.2)
under a concise condition on the live load. In Section 4, a fixed-point iterative method with error estimates is proposed and its applicability and computational efficiency are demonstrated through calculations for two actual bridges; the same calculations quantify both the influence of the nonlinear nonlocal term on the solution and the difference between the two models while explaining several engineering observations. Finally, Section 5 summarizes the main conclusions and proposes several related open problems.

\section{The analytical solution for the ``less stiff'' model (1.5), (1.2) and global properties of the solution}
In this section we study the simplified ``less stiff'' model (1.5), (1.2). After dividing equation (1.5) by $EI$ and denoting the resulting coefficients by $a$ and $c$ and the right-hand side by $\varphi$, we equivalently rewrite problem (1.5),(1.2) as the following linear fourth-order integro-differential problem:
$$
\aligned
&w^{(4)}(x)-aw''(x)+c\int_0^Lw(x)dx=\varphi(x),\ \ \ \ x\in (0,L),\\
&w(0) = w(L) = w''(0) = w''(L) = 0,\\
\endaligned
\eqno (2.1)
$$
where $a,c$ are positive constants, $\varphi\in L^{1}[0,L]$. We first establish the unique solvability of (2.1) and derive an explicit representation of its solution. Although this result was established in [40, Theorem 2.1 and Corollary 2.1], we restate it here for completeness and provide an alternative, more direct proof. We then discuss some global properties of the solution.
%we present here a self-contained statement together with a new proof that is shorter and more direct.

Denote $D(K):=\{w\in W^{4,1}[0,L]:w(0)=w(L)=w''(0)=w''(L)=0\},$ and define the linear operator $K:D(K)\to L^{1}[0,L]$ by
$$
Kw:=w''''-aw''.
$$
To study the Green function of $Kw$, we introduce two second order linear differential operators associated with $K$:
$$
K_{a}w:=-w''+aw,~~~D(K_{a}):=\{w\in W^{2,1}[0,L]:w(0)=w(L)=0\};
$$
$$
K_{0}w:=-w'',~~~D(K_{0}):=\{w\in W^{2,1}[0,L]:w(0)=w(L)=0\}.
$$
The Green functions for $K_{a}$ and $K_{0}$ are
$$G_{a}(t,s)=\left\{
\begin{aligned}
\frac{\sinh (\sqrt a t) \sinh (\sqrt a(L-s))}{\sqrt a\sinh (\sqrt a L)},\ 0\leq t\leq s\leq L,\\
\frac{\sinh (\sqrt as) \sinh (\sqrt a(L-t))}{\sqrt a\sinh (\sqrt a L)},\ 0\leq s\leq t\leq L
\end{aligned}
\right.\eqno (2.2)$$
and
$$G_{0}(t,s)=\left\{ \aligned
\frac{s(L-t)}{L}, \ \ \ \ \ \ \ \  \ 0\leq s\leq t\leq L,\\
\frac{t(L-s)}{L}, \ \ \ \ \ \ \ \  \ 0\leq t\leq s\leq L
\endaligned
\right.\eqno (2.3)
$$ respectively.
Since $K=K_{0}K_{a},$ the Green function of $K$ can be expressed as
$$G(x,s):=\int^{L}_{0}G_{0}(x,t)G_{a}(t,s)dt, \ \ \ \  \ (x,s)\in [0,L]\times [0,L].\eqno (2.4)
$$
As obtained in [40, (2.11)], $G(x,s)$ can be computed as
 $$G(x,s)=\left\{
\begin{aligned}
\frac{x(L-s)}{aL}-\frac{1}{a}\frac{\sinh (\sqrt a x) \sinh (\sqrt a(L-s))}{\sqrt a\sinh (\sqrt a L)},\ 0\leq x\leq s\leq L,\\
\frac{s(L-x)}{aL}-\frac{1}{a}\frac{\sinh (\sqrt a s) \sinh (\sqrt a(L-x))}{\sqrt a\sinh (\sqrt a L)},\ 0\leq s\leq x\leq L.
\end{aligned}
\right.\eqno (2.5)$$
From the nonnegativity of $G_{a}$ and $G_{0}$ it follows that $G(x,s)\geq0$ on $[0,L]\times [0,L]$. Moreover, it is easy to verify from (2.2) that \(G_a(t,s)\) is nonincreasing with respect to the parameter \(a\). Hence, by (2.4) and the nonnegativity of \(G_0\), \(G(x,s)\) is also nonincreasing with
respect to \(a\).

\noindent{\bf Lemma 2.1}\ For any fixed constants $a,c>0$ and any
$\varphi\in L^{1}[0,L]$, the nonlocal boundary value problem (2.1) has a
unique solution
$$
\aligned
w(x)=&\int^{L}_{0}G(x,s)\varphi(s)ds-
\frac{c\int^{L}_{0}G(x,s)ds}
{1+c\int^{L}_{0}\int^{L}_{0}G(t,s)dsdt}
\int^{L}_{0}\int^{L}_{0}G(t,s)\varphi(s)dsdt,
\ \ \ x\in [0,L].
\endaligned\eqno (2.6)
$$

\noindent{\bf Proof.}
Suppose that $w$ is a solution of (2.1), and denote
$$
\psi(x)=\varphi(x)-c\int_0^Lw(t)dt.
$$
Then $w$ satisfies
$$
\aligned
&w^{(4)}(x)-aw''(x)=\psi(x),\ \ \ \ x\in (0,L),\\
&w(0) = w(L) = w''(0) = w''(L) = 0.
\endaligned
\eqno (2.7)
$$
By (2.4) and (2.5), the solution of (2.7) can be expressed as follows:
$$
\aligned
w(x)&=\int^{L}_{0}G(x,s)\psi(s)ds\\
&=\int^{L}_{0}G(x,s)\varphi(s)ds-
c\int^{L}_{0}G(x,s)ds\int^{L}_{0}w(t)dt,
\ \ \ x\in [0,L].
\endaligned\eqno (2.8)
$$
Integrating both sides of (2.8) from $0$ to $L$, and noting that
$c>0$ and $G(t,s)\geq0$, we obtain
$$
\aligned
\int^{L}_{0}w(t)dt
=\frac{\int^{L}_{0}\int^{L}_{0}G(t,s)\varphi(s)dsdt}
{1+c\int^{L}_{0}\int^{L}_{0}G(t,s)dsdt}.
\endaligned\eqno (2.9)
$$
Substituting (2.9) into (2.8) yields
$$
\aligned
w(x)=&\int^{L}_{0}G(x,s)\varphi(s)ds-
\frac{c\int^{L}_{0}G(x,s)ds}
{1+c\int^{L}_{0}\int^{L}_{0}G(t,s)dsdt}
\int^{L}_{0}\int^{L}_{0}G(t,s)\varphi(s)dsdt,
\ \ \ x\in [0,L],
\endaligned
$$
so every solution of (2.1) must be given by (2.6).

Conversely, define $w$ by (2.6), and denote the right-hand side of
(2.9) by $I_{\varphi}$. By the defining properties of the Green function
$G$, we have $w\in D(K)$ and
$$
w^{(4)}(x)-aw''(x)=\varphi(x)-cI_{\varphi}.
$$
Moreover, direct integration of (2.6) gives
$$
\int^{L}_{0}w(t)dt=I_{\varphi}.
$$
Consequently,
$$
w^{(4)}(x)-aw''(x)+c\int^{L}_{0}w(t)dt
=\varphi(x),
$$
and hence $w$ is a solution of (2.1). Together with the preceding
uniqueness argument, this completes the proof. \hfill{$\Box$}

\noindent{\bf Remark 2.1} As demonstrated in [40], the explicit solution (2.6) provides not only a crucial theoretical understanding of the fundamental behavior of the ``less stiff'' model (1.5), (1.2) but also an exact benchmark, serving as a definitive criterion for evaluating the numerical methods developed in [28-33]. The proof presented here is more concise.

\noindent{\bf Lemma 2.2}\
For any $\varphi\in L^{1}[0,L]$, the solution $w$ of problem (2.1) given
by (2.6) satisfies
$$
|\int_0^L w(x)dx|
\leq
\frac{L^2}{8a \left[1 + c \left( \frac{L^3}{12a}
- \frac{L}{a^2}
+\frac{2}{a^{5/2}} \tanh\left(\frac{\sqrt{a}L}{2}\right)
\right) \right]}\|\varphi\|_{1}
:=E(a)\|\varphi\|_{1}.
\eqno (2.10)
$$
Moreover, \(E(a)\) is nonincreasing with respect to \(a\).

\noindent{\bf Proof.} Based upon (2.6), by explicit calculation and
interchanging the order of integration we have
$$
\aligned
|\int_0^L w(x)dx|
&=|\int_0^L[\int^{L}_{0}G(x,s)\varphi(s)ds-
\frac{c\int^{L}_{0}G(x,s)ds}{1+c\int^{L}_{0}\int^{L}_{0}G(x,s)dsdx}
\int^{L}_{0}\int^{L}_{0}G(x,s)\varphi(s)dsdx]dx|      \\
&=|\int_0^L\int^{L}_{0}G(x,s)\varphi(s)dsdx-
\frac{c\int^{L}_{0}\int^{L}_{0}G(x,s)dsdx}
{1+c\int^{L}_{0}\int^{L}_{0}G(x,s)dsdx}
\int^{L}_{0}\int^{L}_{0}G(x,s)\varphi(s)dsdx|      \\
&=\frac{|\int_0^L\int^{L}_{0}G(x,s)\varphi(s)dsdx|}
{1+c\int^{L}_{0}\int^{L}_{0}G(x,s)dsdx}
=\frac{|\int_0^L(\int^{L}_{0}G(x,s)dx)\varphi(s)ds|}
{1+c\int^{L}_{0}\int^{L}_{0}G(x,s)dsdx}\\
&\leq
\frac{\int_0^L(\int^{L}_{0}G(x,s)dx)|\varphi(s)|ds}
{1+c\int^{L}_{0}\int^{L}_{0}G(x,s)dsdx}
\leq
\frac{\underset{s\in [0,L]}\max\int^{L}_{0}G(x,s)dx\|\varphi\|_{1}}
{1+c\int^{L}_{0}\int^{L}_{0}G(x,s)dsdx}.
\endaligned\eqno (2.11)
$$
Since
$$
\underset{s\in [0,L]}\max\int^{L}_{0}G(x,s)dx
=
\frac{L^2}{8a} - \frac{1}{a^2}
+ \frac{1}{a^2}\operatorname{sech}\left(\frac{\sqrt{a}L}{2}\right)
$$
and
$$
\int^{L}_{0}\int^{L}_{0}G(x,s)dsdx
=
\frac{L^3}{12a} - \frac{L}{a^2}
+ \frac{2}{a^{5/2}} \tanh\left(\frac{\sqrt{a}L}{2}\right),
$$
then by (2.11) and the nonnegativity of
$1-\operatorname{sech}(x)$ we have
$$
\aligned
|\int_0^L w(x)dx|
&\leq
\frac{\frac{L^2}{8a} - \frac{1}{a^2}
+ \frac{1}{a^2}\operatorname{sech}\left(\frac{\sqrt{a}L}{2}\right)}
{1+c\left( \frac{L^3}{12a} - \frac{L}{a^2}
+\frac{2}{a^{5/2}} \tanh\left(\frac{\sqrt{a}L}{2}\right) \right)}
\|\varphi\|_{1}\\
&\leq
\frac{L^2}{8a \left[1 + c \left( \frac{L^3}{12a}
- \frac{L}{a^2}
+\frac{2}{a^{5/2}} \tanh\left(\frac{\sqrt{a}L}{2}\right)
\right) \right]}\|\varphi\|_{1}
=E(a)\|\varphi\|_{1},
\endaligned
$$
thus, (2.10) is proved. To prove that \(E(a)\) is nonincreasing with
respect to \(a\), it is enough to show that
$$
a\left[1+c\left(
\frac{L^3}{12a}
-\frac{L}{a^2}
+\frac{2}{a^{5/2}}\tanh\left(\frac{\sqrt aL}{2}\right)
\right)\right]
\eqno (2.12)
$$
is nondecreasing with respect to \(a\). By direct differentiation, we have
$$
\begin{aligned}
&\frac{d}{da}\left[
a\left(\frac{L^3}{12a}-\frac{L}{a^2}+\frac{2}{a^{5/2}}\tanh\left(\frac{\sqrt aL}{2}\right)\right)\right]=
\frac{2\frac{\sqrt aL}{2}+\frac{\sqrt aL}{2}\operatorname{sech}^2\!\left(\frac{\sqrt aL}{2}\right)-3\tanh\!\left(\frac{\sqrt aL}{2}\right)}{a^{5/2}} \end{aligned}
\eqno (2.13)
$$
Since it is easy to verify that
\[
2x+x\operatorname{sech}^2x-3\tanh x\geq0,\qquad x\geq0,
\]
taking \(x=\frac{\sqrt aL}{2}\), we obtain that the derivative in
(2.13) is nonnegative. Hence the expression in (2.12) is
nondecreasing with respect to \(a\). Therefore \(E(a)\) is nonincreasing
with respect to \(a\) and the proof of the lemma is complete.\hfill{$\Box$}

\noindent{\bf Remark 2.2} While Lemma 2.2 provides an a priori estimate
for the integral of the solution, (2.9) further shows
that, when $\varphi\geq0$, the solution $w$ of problem (2.1) satisfies
$$
\int_0^L w(x)dx
=\frac{\int_0^L\int^{L}_{0}G(x,s)\varphi(s)dsdx}
{1+c\int^{L}_{0}\int^{L}_{0}G(x,s)dsdx}
\geq0
$$
without any restriction on the parameters. In contrast, the pointwise
positivity of $w$, i.e., $w(x)\geq0$ for all $x\in[0,L]$, under the same
condition $\varphi\geq0,$ requires a certain condition on the parameters;
such a restrictive sufficient condition has been established in
[40, Theorem 3.1].

\noindent{\bf Lemma 2.3}\
For any $\varphi\in L^{1}[0,L]$, the solution $w$ of problem (2.1) given
by (2.6) satisfies
$$
w''(x)
=
-\int^{L}_{0}G_{a}(x,s)\varphi(s)ds+
\frac{c\int^{L}_{0}\int^{L}_{0}G(x,s)\varphi(s)dsdx}
{1+c\int^{L}_{0}\int^{L}_{0}G(x,s)dsdx}
\int^{L}_{0}G_{a}(x,s)ds,
\eqno (2.14)
$$
and
$$
\int_0^L |w''(x)|dx=\|w''\|_{1}
\leq
\frac{
1+
c\left(
\frac{5L^3}{24a}
-\frac{L}{a^2}
+
\frac{1}{a^{5/2}}
\left(
2-\frac{aL^2}{4}
\right)
\tanh\!\left(\frac{\sqrt aL}{2}\right)
\right)
}
{
a\left[1+c\left(\frac{L^3}{12a}
-\frac{L}{a^2}
+\frac{2}{a^{5/2}}\tanh\!\left(\frac{\sqrt aL}{2}\right)\right)\right]}
\|\varphi\|_{1}
:=D(a)\|\varphi\|_{1}.
\eqno (2.15)
$$
Moreover, \(D(a)\) is nonincreasing with respect to \(a\).

\noindent{\bf Proof.} Based upon (2.6) and (2.4), by explicit calculation
we can obtain (2.14). Then, by (2.14) and (2.2) and by interchanging the
order of integration, we have
$$
\aligned
\int_0^L |w''(x)|dx
&=\int^{L}_{0}\left|-\int^{L}_{0}G_{a}(x,s)\varphi(s)ds+
\frac{c\int^{L}_{0}\int^{L}_{0}G(x,s)\varphi(s)dsdx}
{1+c\int^{L}_{0}\int^{L}_{0}G(x,s)dsdx}
\int^{L}_{0}G_{a}(x,s)ds\right|dx      \\
&\leq
\int_0^L\int^{L}_{0}G_{a}(x,s)|\varphi(s)|dsdx+
\frac{c\int^{L}_{0}\int^{L}_{0}G(x,s)|\varphi(s)|dsdx}
{1+c\int^{L}_{0}\int^{L}_{0}G(x,s)dsdx}
\int^{L}_{0}\int^{L}_{0}G_{a}(x,s)dsdx     \\
&=
\int_0^L(\int^{L}_{0}G_{a}(x,s)dx)|\varphi(s)|ds+
\frac{c\int^{L}_{0}(\int^{L}_{0}G(x,s)dx)|\varphi(s)|ds}
{1+c\int^{L}_{0}\int^{L}_{0}G(x,s)dsdx}
\int^{L}_{0}\int^{L}_{0}G_{a}(x,s)dsdx     \\
&\leq
\underset{s\in [0,L]}\max\int^{L}_{0}G_{a}(x,s)dx\|\varphi\|_{1}
+
\frac{\underset{s\in [0,L]}\max\int^{L}_{0}G(x,s)dx\|\varphi\|_{1}}
{1+c\int^{L}_{0}\int^{L}_{0}G(x,s)dsdx}
c\int^{L}_{0}\int^{L}_{0}G_{a}(x,s)dsdx.
\endaligned\eqno (2.16)
$$
Since
$$
\underset{s\in [0,L]}\max\int^{L}_{0}G_a(x,s)dx
=
\frac{1}{a}\left(1 - \operatorname{sech}\left(\frac{\sqrt{a}L}{2}\right)\right)
$$
and
$$
\int^{L}_{0}\int^{L}_{0}G_a(x,s)dsdx
=
\frac{L}{a} - \frac{2}{a\sqrt{a}}
\tanh\left(\frac{\sqrt{a}L}{2}\right),
$$
by (2.16) together with Lemma 2.2 and the nonnegativity of
$\operatorname{sech}(x)$ we have
$$
\aligned
\int_0^L |w''(x)|dx
&\leq
\left[
\frac{1}{a}\left(1 - \operatorname{sech}\left(\frac{\sqrt{a}L}{2}\right)\right)
+
\frac{L^2c\left(\frac{L}{a}
-\frac{2}{a\sqrt{a}}\tanh\left(\frac{\sqrt{a}L}{2}\right)\right)}
{8a \left[1 + c \left( \frac{L^3}{12a}
- \frac{L}{a^2}
+\frac{2}{a^{5/2}} \tanh\left(\frac{\sqrt{a}L}{2}\right)
\right) \right]}
\right]\|\varphi\|_{1}\\
&\leq
\frac{
1+
c\left(
\frac{5L^3}{24a}
-\frac{L}{a^2}
+
\frac{1}{a^{5/2}}
\left(
2-\frac{aL^2}{4}
\right)
\tanh\!\left(\frac{\sqrt aL}{2}\right)
\right)
}
{
a\left[1+c\left(\frac{L^3}{12a}
-\frac{L}{a^2}
+\frac{2}{a^{5/2}}\tanh\!\left(\frac{\sqrt aL}{2}\right)\right)\right]}
\|\varphi\|_{1}
=D(a)\|\varphi\|_{1},
\endaligned\eqno (2.17)
$$
thus, (2.15) is proved. To prove that \(D(a)\) is nonincreasing with
respect to \(a\), we first notice that its denominator is exactly
(2.12) and its nondecreasing property has been proved in Lemma 2.2.
Therefore, it remains only to show that
$$
1+c\left(
\frac{5L^3}{24a}
-\frac{L}{a^2}
+\frac{1}{a^{5/2}}
\left(2-\frac{aL^2}{4}\right)
\tanh\left(\frac{\sqrt aL}{2}\right)
\right)
\eqno (2.18)
$$
is nonincreasing with respect to \(a\). Indeed, according to the explicit
computations in Lemma 2.2 and in the proof above, (2.18) can be written as
$$
1+c\left[
\int_0^L\int_0^L G(x,s)\,dsdx
+
\frac{L^2}{8}
\int_0^L\int_0^L G_a(x,s)\,dsdx
\right].
\eqno (2.19)
$$
Since both \(G_a(x,s)\) and \(G(x,s)\) are nonincreasing with respect to
\(a\), (2.19), equivalently (2.18), is nonincreasing with respect to
\(a\). Therefore \(D(a)\) is nonincreasing with respect to \(a\), and the
proof of the lemma is complete.\hfill{$\Box$}

\noindent{\bf Remark 2.3} By (2.14), the representation of $w''$ in Lemma 2.3 can be rewritten, using (2.9), as
$$
w''(x)=-\int^{L}_{0}G_{a}(x,s)\varphi(s)ds
+c\int^{L}_{0}w(t)dt\int^{L}_{0}G_{a}(x,s)ds.
$$
When $\varphi\geq0$, Remark 2.2 shows that
$\int^{L}_{0}w(t)dt\geq0$. Hence, the first term on the right-hand side
is nonpositive whereas the second is nonnegative, so
$\varphi\geq0$ alone does not ensure that $w''\leq0$. A restrictive
sufficient condition on the parameters ensuring that
$w''(x)\leq0$ for all $x\in[0,L]$ was established in
[40, Theorem 3.1]. Under this condition, the deflection curve predicted
by the simplified ``less stiff'' model (1.5),(1.2) has no inflection
points. 

\section{Existence and uniqueness of solutions to the classical Melan equation (1.4) with (1.2)}

In this section, we study the original nonlinear nonlocal problem (1.4),(1.2). After dividing equation (1.4) by \(EI\) and introducing the parameters below, we write it equivalently as the following fourth-order equation subject to the simply supported boundary conditions (1.2):
$$
w''''(x)-(a+b\int_0^L w(x)dx)w''(x)+c\int_0^L w(x)dx=p(x),\ \ \ \ x\in (0,L),\eqno (3.1)
$$
where $a=\frac{H}{EI},b=\frac{1}{EI}\frac{E_{c}A_{c}}{L_{c}}\frac{q}{H},c=\frac{1}{EI}\frac{q^{2}}{H^{2}}\frac{E_{c}A_{c}}{L_{c}}$ are positive parameters, $p=\frac{P(x)}{EI}\in L^{1}[0,L]$.
%$$Gw=w^{(4)}-Mw''+N\int_0^Lw(t)dt  \eqno (3.2)
%$$ on $D(G):=\{w\in W^{4,1}[0,L]:w(0)=w(L)=w''(0)=w''(L)=0\}.$
For problem (3.1),(1.2), we first derive a priori bounds for any solution with nonnegative integral and its corresponding nonlinear term, using the estimates established for the associated linear problem in Section 2.

\noindent{\bf Lemma 3.1}\
Let \(w\in W^{4,1}[0,L]\) be a solution of the BVP (3.1),(1.2) for a
given \(p\in L^1[0,L]\) satisfying
$$
\frac{2a^2 \left[1+c\left(\frac{L^3}{12a}-\frac{L}{a^2}
+\frac{2}{a^{5/2}}\tanh\!\left(\frac{\sqrt aL}{2}\right)\right)\right]^2}
{ bL^2\left[ 1+ c\left( \frac{5L^3}{24a} -\frac{L}{a^2}
+ \frac{1}{a^{5/2}} \left( 2-\frac{aL^2}{4} \right)
  \tanh\!\left(\frac{\sqrt aL}{2}\right) \right) \right]}
>
\|p\|_{1}.
\eqno (H)
$$
If \(\int_0^L w(x)dx\geq0\), then
$$
\left|\int_0^L w(x)dx\right|
\leq E(a)\|p\|_1,
\qquad
\|w''\|_1\leq D(a)\|p\|_1,
$$
where $E(a)$ and $D(a)$ are defined in (2.10) and (2.15),
respectively.
Moreover, for
$$
\varphi=p+b\left(\int_0^Lw(x)dx\right)w'',
$$
we have
$$
\|\varphi\|_1\leq \frac54\|p\|_1.
\eqno (3.2)
$$

\noindent{\bf Proof.}
Let
$$
I=\int_0^L w(x)dx,
$$
then problem (3.1),(1.2) can be rewritten as the following linear nonlocal problem:
$$
w''''(x)-(a+bI)w''(x)+c\int_0^Lw(x)dx=p(x),
\qquad x\in(0,L),
$$
with the boundary conditions (1.2), which has been studied in Section 2.
Since \(I\geq0\), we have \(a+bI\geq a\). Applying Lemma 2.2 and Lemma
2.3 with \(a\) replaced by \(a+bI\), and using the fact that \(E(a)\) and
\(D(a)\) are nonincreasing with respect to \(a\), we obtain
$$
\left|\int_0^L w(x)dx\right|
\leq E(a+bI)\|p\|_1
\leq E(a)\|p\|_1,
$$
and
$$
\|w''\|_1
\leq D(a+bI)\|p\|_1
\leq D(a)\|p\|_1 .
$$
By the definitions of \(E(a)\) and \(D(a)\), hypothesis (H) is equivalent
to
$$
4bE(a)D(a)\|p\|_1<1,\eqno (3.3)
$$
therefore,
$$
b\left|\int_0^L w(x)dx\right|\|w''\|_1
\leq bE(a)D(a)\|p\|_1^2
\leq
\frac14\|p\|_1 .
$$
Consequently,
$$
\begin{aligned}
\|\varphi\|_1
&=
\left\|p+b\left(\int_0^Lw(x)dx\right)w''\right\|_1\leq
\|p\|_1+
b\left|\int_0^Lw(x)dx\right|\|w''\|_1\leq
\frac54\|p\|_1 .
\end{aligned}
$$
The proof is complete.\hfill{$\Box$}

Motivated by the expression estimated in Lemma 3.1, we now formulate
nonlinear BVP (3.1),(1.2) as a fixed point problem. For each \(\varphi\in L^1[0,L]\), let \(w\) be the unique solution of the linear nonlocal problem (2.1), which can be expressed explicitly by (2.6). Define
\(A:L^1[0,L]\to L^1[0,L]\) by
$$
(A\varphi)(x)
=
p(x)+b\left(\int_0^L w(t)dt\right)w''(x),
\eqno (3.4)
$$
where \(w''\) is given by (2.14), again in terms of \(\varphi\).

We then have the following equivalence between the BVP (3.1),(1.2) and
the fixed point equation \(A\varphi=\varphi\), which follows directly from
the definition of \(A\) and the representation formula (2.6).

\noindent{\bf Lemma 3.2}\ If the BVP (3.1),(1.2) has a solution $w(x)\in W^{4,1}[0,L]$ for a given $p\in L^{1}[0,L]$, then the function
$$
\varphi(x)=p(x)+b\int_0^L w(x)dxw''(x)\eqno (3.5)
$$ is a fixed point of the operator $A$ defined by (3.4). Conversely, if the function $\varphi$ is a fixed point of the operator $A$, then the function $w(x)\in W^{4,1}[0,L]$ determined from (2.6) in terms of $\varphi$ is a solution of the BVP (3.1),(1.2).

By Lemma 3.2, to prove the existence of a solution of the BVP
(3.1),(1.2), it is enough to show that the operator \(A\) has a fixed
point. We shall do this by applying the Banach fixed point theorem on a
closed ball in \(L^1[0,L]\). Let \(B[0,2\|p\|_1]\) denote the closed ball in \(L^1[0,L]\) centered at
\(0\) with radius \(2\|p\|_1\), that is,
$$
B[0,2\|p\|_1]=\{\varphi\in L^1[0,L]\;:\;\|\varphi\|_1\leq2\|p\|_1\}.
$$
The first main result is as follows:

\noindent{\bf Theorem 3.1}\
For any fixed constants \(a,b,c>0\) and \(p\in L^{1}[0,L]\), assume that
hypothesis (H) holds. Then the operator \(A\) has a unique fixed point
\(\varphi\in B[0,2\|p\|_1]\). Consequently, the BVP (3.1),(1.2) has a
corresponding solution \(w\in W^{4,1}[0,L]\hookrightarrow\hookrightarrow
C^{2}[0,L]\), determined by (2.6) in terms of \(\varphi\). Moreover, this
solution satisfies
$$
\|w\|_{1}
\leq
\frac{1+2c\left[\frac{L^3}{12a} - \frac{L}{a^2}
+ \frac{2}{a^{5/2}} \tanh\left(\frac{\sqrt{a}L}{2}\right)\right]}
{1+c\left[\frac{L^3}{12a} - \frac{L}{a^2}
+ \frac{2}{a^{5/2}} \tanh\left(\frac{\sqrt{a}L}{2}\right)\right]}
\frac{L^2}{4a}\|p\|_{1}
$$
and
$$
\|w''\|_{1}
\leq
\frac{
1+c\left(\frac{5L^3}{24a}-\frac{L}{a^2}
+\frac{1}{a^{5/2}}
\left(2-\frac{aL^2}{4}\right)
\tanh\!\left(\frac{\sqrt aL}{2}\right)\right)}
{a\left[1+c\left(\frac{L^3}{12a}-\frac{L}{a^2}
+\frac{2}{a^{5/2}}\tanh\!\left(\frac{\sqrt aL}{2}\right)\right)\right]}
2\|p\|_{1}.
$$

\noindent{\bf Proof.}
We show that \(A\) is a contraction mapping on \(B[0,2\|p\|_1]\). Take
any \(\varphi\in B[0,2\|p\|_1]\), then the linear nonlocal problem (2.1) admits a unique solution $w$, whose explicit form is given by (2.6), and $w''$ is given by (2.14).  Hence, by (3.4) and (2.10) in Lemma 2.2 and (2.15) in Lemma 2.3 we have
$$
\aligned
\|A\varphi\|_1
&=
\left\|p+b\left(\int_0^Lw(t)dt\right)w''\right\|_1\leq
\|p\|_1+
b\left|\int_0^Lw(t)dt\right|\|w''\|_1\\
&\leq
\|p\|_1+bE(a)D(a)\|\varphi\|_1^2.
\endaligned\eqno (3.6)
$$
Since \(\varphi \in B[0,2\|p\|_{1}]\), that is
\(\|\varphi\|_{1}\leq 2\|p\|_{1}\), combining this with the equivalent form
(3.3) of hypothesis (H), it follows from (3.6) that
$$
\|A\varphi\|_1\leq\|p\|_1+4bE(a)D(a)\|p\|_1^2\leq2\|p\|_1,
$$
which means that $A$ maps $B[0,2\|p\|_{1}]$ into itself.

Next, take \(\varphi_1,\varphi_2\in B[0,2\|p\|_1]\), and let \(w_1,w_2\)
be the corresponding solutions of the linear nonlocal problem (2.1). Then, by (3.4) and (2.10) in Lemma 2.2 and (2.15) in Lemma 2.3, we have
$$
\aligned
&\|A\varphi_{2}-A\varphi_{1}\|_{1}=\int^{L}_{0}|A\varphi_{2}-A\varphi_{1}|dx\\
&= \int^{L}_{0}|[p(x)+b\int_0^L w_{2}(x)dxw_{2}''(x)]-[p(x)+b\int_0^L w_{1}(x)dxw_{1}''(x)]|dx\\
&= \int^{L}_{0}|b\int_0^L w_{2}(x)dxw_{2}''(x)-b\int_0^L w_{1}(x)dxw_{2}''(x)+b\int_0^L w_{1}(x)dxw_{2}''(x)-b\int_0^L w_{1}(x)dxw_{1}''(x)|dx\\
&= \int^{L}_{0}|b\int_0^L (w_{2}(x)-w_{1}(x))dxw_{2}''(x)+b\int_0^L w_{1}(x)dx(w_{2}''(x)-w_{1}''(x))|dx\\
&\leq b|\int_0^L (w_{2}(x)-w_{1}(x))dx|\int^{L}_{0}|w_{2}''(x)|dx+b|\int_0^L w_{1}(x)dx|\int_0^L|w_{2}''(x)-w_{1}''(x)|dx\\
&\leq
bE(a)\|\varphi_{2}-\varphi_{1}\|_{1}D(a)\|\varphi_{2}\|_{1}
+
bE(a)\|\varphi_{1}\|_{1}D(a)\|\varphi_{2}-\varphi_{1}\|_{1}\\
&\leq
4bE(a)D(a)\|p\|_{1}\|\varphi_{2}-\varphi_{1}\|_{1}\\
&=
\|p\|_{1}
\frac{
bL^2\left[
1+
c\left(
\frac{5L^3}{24a}
-\frac{L}{a^2}
+
\frac{1}{a^{5/2}}
\left(
2-\frac{aL^2}{4}
\right)
\tanh\!\left(\frac{\sqrt aL}{2}\right)
\right)
\right]
}{
2a^2
\left[
1+c\left(
\frac{L^3}{12a}
-\frac{L}{a^2}
+\frac{2}{a^{5/2}}
\tanh\!\left(\frac{\sqrt aL}{2}\right)
\right)
\right]^2
}
\|\varphi_{2}-\varphi_{1}\|_{1},\\
\endaligned\eqno (3.7)
  $$
Combining (3.7) with the equivalent form (3.3) of hypothesis (H), we have
$$\|A\varphi_{2}-A\varphi_{1}\|_{1}\leq q\|\varphi_{2}-\varphi_{1}\|_{1},\eqno (3.8)$$
where $$q=4bE(a)D(a)\|p\|_1=\|p\|_{1}\frac{
bL^2\left[
1+
c\left(
\frac{5L^3}{24a}
-\frac{L}{a^2}
+
\frac{1}{a^{5/2}}
\left(
2-\frac{aL^2}{4}
\right)
\tanh\!\left(\frac{\sqrt aL}{2}\right)
\right)
\right]
}{
2a^2
\left[1+c\left(\frac{L^3}{12a}-\frac{L}{a^2}+\frac{2}{a^{5/2}}\tanh\!\left(\frac{\sqrt aL}{2}\right)\right)\right]^2
}<1,\eqno (3.9)
 $$
therefore \(A\) is a contraction mapping on \(B[0,2\|p\|_1]\). By the Banach
fixed point theorem, \(A\) has a unique fixed point
\(\varphi\in B[0,2\|p\|_1]\), and then by Lemma 3.2 the function \(w\) determined by (2.6) in terms of \(\varphi\) is a solution of the BVP (3.1),(1.2).

It remains to prove the estimates for \(w\). By a direct computation we have:
$$
\aligned
&\|w\|_{1}=\int^{L}_{0}|\int^{L}_{0}G(x,s)\varphi(s)ds-
\frac{c\int^{L}_{0}G(x,s)ds}{1+c\int^{L}_{0}\int^{L}_{0}G(x,s)dsdx}\int^{L}_{0}\int^{L}_{0}G(x,s)\varphi(s)dsdx|dx\\
&  \leq \int_0^L\int^{L}_{0}G(x,s)|\varphi(s)|dsdx+
\frac{c\int^{L}_{0}\int^{L}_{0}G(x,s)dsdx}{1+c\int^{L}_{0}\int^{L}_{0}G(x,s)dsdx}\int^{L}_{0}\int^{L}_{0}G(x,s)|\varphi(s)|dsdx\\
& =\frac{1+2c\int^{L}_{0}\int^{L}_{0}G(x,s)dsdx}{1+c\int^{L}_{0}\int^{L}_{0}G(x,s)dsdx}\int^{L}_{0}(\int^{L}_{0}G(x,s)dx)|\varphi(s)|ds\\ %=\frac{1+2c\int^{L}_{0}\int^{L}_{0}G(x,s)dsdx}{1+c\int^{L}_{0}\int^{L}_{0}G(x,s)dsdx}\int^{L}_{0}(\int^{L}_{0}G(x,s)dx)|\varphi(s)|ds\\
&
\leq\frac{1+2c\int^{L}_{0}\int^{L}_{0}G(x,s)dsdx}{1+c\int^{L}_{0}\int^{L}_{0}G(x,s)dsdx}\underset{s\in [0,L]}\max\int^{L}_{0}G(x,s)dx\|\varphi\|_{1}\\
%&\ \ \ \ \ \ \
%= \frac{1+2c[\frac{L^3}{12a} - \frac{L}{a^2} + \frac{2}{a^{5/2}} \tanh\left(\frac{\sqrt{a}L}{2}\right)]}{1+c[\frac{L^3}{12a} - \frac{L}{a^2} + \frac{2}{a^{5/2}} \tanh\left(\frac{\sqrt{a}L}{2}\right)]}(\frac{L^2}{8a} - \frac{1}{a^2} + \frac{1}{a^2}\operatorname{sech}\left(\frac{\sqrt{a}L}{2}\right))\|\varphi\|_{1}\\
&
\leq \frac{1+2c[\frac{L^3}{12a} - \frac{L}{a^2} + \frac{2}{a^{5/2}} \tanh\left(\frac{\sqrt{a}L}{2}\right)]}{1+c[\frac{L^3}{12a} - \frac{L}{a^2} + \frac{2}{a^{5/2}} \tanh\left(\frac{\sqrt{a}L}{2}\right)]}\frac{L^2}{8a}\|\varphi\|_{1}
%&\ \ \ \ \ \ \
%\leq \frac{1+2c[\frac{L^3}{12a} - \frac{L}{a^2} + \frac{2}{a^{5/2}} \tanh\left(\frac{\sqrt{a}L}{2}\right)]}{1+c[\frac{L^3}{12a} - \frac{L}{a^2} + \frac{2}{a^{5/2}} \tanh\left(\frac{\sqrt{a}L}{2}\right)]}\frac{L^2}{8a}2\|p\|_{1}\\
\leq \frac{1+2c[\frac{L^3}{12a} - \frac{L}{a^2} + \frac{2}{a^{5/2}} \tanh\left(\frac{\sqrt{a}L}{2}\right)]}{1+c[\frac{L^3}{12a} - \frac{L}{a^2} + \frac{2}{a^{5/2}} \tanh\left(\frac{\sqrt{a}L}{2}\right)]}\frac{L^2}{4a}\|p\|_{1};\\
\endaligned
\eqno (3.10)
$$
Moreover, since \(\varphi\in B[0,2\|p\|_1]\), we have
\(\|\varphi\|_1\leq 2\|p\|_1\). Hence, by (2.15) in Lemma 2.3,
$$
\|w''\|_1\leq
D(a)\|\varphi\|_1
\leq
2D(a)\|p\|_1,
$$
which is exactly the desired estimate for \(\|w''\|_1\). This completes
the proof of the theorem.\hfill{$\Box$}

Although Theorem 3.1 does not directly give global uniqueness for the
nonlinear BVP (3.1),(1.2), its fixed point argument yields the following restricted uniqueness result:

\noindent{\bf Theorem 3.2}\
Assume that hypothesis (H) holds, and let \(w\) be the solution obtained
in Theorem 3.1. If \(\widetilde w\in W^{4,1}[0,L]\) is any solution of the
BVP (3.1),(1.2) satisfying
\[
\int_0^L\widetilde w(x)dx\geq0,
\]
then
\[
\widetilde w=w,
\]
which implies that a solution of the BVP (3.1),(1.2) in the class
$$
\mathcal X_+
=
\left\{
u\in W^{4,1}[0,L]\;:\;
\int_0^L u(x)dx\geq0
\right\}
$$
is unique whenever it exists.

\noindent{\bf Proof.}
Let \(\widetilde w\in W^{4,1}[0,L]\) be a solution of (3.1),(1.2)
satisfying \(\int_0^L\widetilde w(x)dx\geq0\). By Lemma 3.1, the function
$$
\widetilde\varphi
=
p+b\left(\int_0^L\widetilde w(x)dx\right)\widetilde w''
$$
satisfies
$$
\|\widetilde\varphi\|_1\leq\frac54\|p\|_1\leq2\|p\|_1.
$$
Hence \(\widetilde\varphi\in B[0,2\|p\|_1]\). By Lemma 3.2,
\(\widetilde\varphi\) is a fixed point of \(A\). Since \(A\) has a unique
fixed point in \(B[0,2\|p\|_1]\), we have
\(\widetilde\varphi=\varphi\), where \(\varphi\) is the fixed point
obtained in Theorem 3.1. Therefore, by (2.6), \(\widetilde w=w\). As \(\widetilde w\) was chosen as an arbitrary solution of (3.1),(1.2)
in \(\mathcal X_+\), the uniqueness in this class follows. The
proof is complete.\hfill{$\Box$}

%Corollary 3.1 shows that any solution satisfying the nonnegative-integral condition is exactly the solution corresponding to the unique fixed point obtained in Theorem 3.1. We next show that, if the live load is restricted to be nonnegative, then this solution indeed satisfies the nonnegative-integral condition.
\noindent{\bf Remark 3.1}\ Theorem 3.2 establishes global uniqueness in \(\mathcal X_+\) whenever a solution in this class exists, without imposing any smallness restriction on the solution. From the physical and engineering point of view, the restriction defining \(\mathcal X_+\) is natural and makes this result the relevant global uniqueness statement for the deflection curve. Indeed, usual live loads are mainly downward traffic loads, so the corresponding deflection is expected to have nonnegative total displacement. It should be emphasized that Theorem 3.2 itself does not require $p\geq0$ and may therefore also apply to certain sign-changing loads involving small upward effects such as wind lift, provided that the resulting solution remains in \(\mathcal X_+\).

Theorem 3.2 shows that any solution in \(\mathcal X_+\) is exactly the
solution corresponding to the unique fixed point obtained in Theorem 3.1.
We next prove that, for nonnegative live loads, this fixed-point solution actually belongs to \(\mathcal X_+\).

\noindent{\bf Theorem 3.3}\
Assume that hypothesis (H) holds and \(p\geq0\). Then the solution \(w\)
obtained in Theorem 3.1 satisfies
$$
\int_0^L w(x)dx\geq0,
$$
that is,  \(w\in\mathcal X_+\). Consequently, by Theorem 3.2, the BVP (3.1),(1.2) has a unique solution with nonnegative integral.

\noindent{\bf Proof.}
Let \(\varphi\in B[0,2\|p\|_1]\) be the fixed point obtained in Theorem
3.1, and let \(w\) be the corresponding solution of the BVP
(3.1),(1.2). Since \(\varphi\) is a fixed point of \(A\), we have
$$
\varphi(s)=p(s)+b\left(\int_0^Lw(x)dx\right)w''(s).
\eqno (3.11)
$$
Using the relation between \(w\) and \(\varphi\) given by (2.6) and
integrating both sides over \([0,L]\), we obtain
$$
\int_0^Lw(x)dx
=
\frac{1}{1+c\int_0^L\int_0^LG(x,s)dsdx}
\int_0^L\int_0^LG(x,s)\varphi(s)dsdx,
$$
and hence
$$
\left(\int_0^Lw(x)dx\right)
\left[
1+c\int_0^L\int_0^LG(x,s)dsdx
\right]
=
\int_0^L\int_0^LG(x,s)\varphi(s)dsdx .
\eqno (3.12)
$$
Substituting (3.11) into the right-hand side of (3.12), we obtain
$$
\begin{aligned}
&\left(\int_0^Lw(x)dx\right)
\left[
1+c\int_0^L\int_0^LG(x,s)dsdx
-
b\int_0^L\int_0^LG(x,s)w''(s)dsdx
\right] \\
&\quad =
\int_0^L\int_0^LG(x,s)p(s)dsdx .
\end{aligned}
\eqno (3.13)
$$
Since \(G(x,s)\geq0\) and \(p\geq0\), the right-hand side of (3.13) is
nonnegative. Therefore, to prove
$$
\int_0^Lw(x)dx\geq0,
$$
it is enough to show that
$$
1+c\int_0^L\int_0^LG(x,s)dsdx
-
b\int_0^L\int_0^LG(x,s)w''(s)dsdx>0.\eqno (3.14)
$$
Indeed,
$$
\begin{aligned}
&b\int_0^L\int_0^LG(x,s)w''(s)dsdx\leq
b\int_0^L\left(\int_0^LG(x,s)dx\right)|w''(s)|ds\\
&\quad\quad\quad\quad\quad\quad\quad\quad\quad\quad\quad\quad
\leq b\max_{s\in[0,L]}\int_0^LG(x,s)dx\, \|w''\|_1.
\end{aligned}\eqno (3.15)
$$
By (2.11) and the definition of \(E(a)\) in Lemma 2.2, we have
$$
\max_{s\in[0,L]}\int_0^LG(x,s)dx
\leq
E(a)\left[
1+c\int_0^L\int_0^LG(x,s)dsdx
\right].
$$
Combining this with the estimate
$\|w''\|_1\leq2D(a)\|p\|_1$ obtained in Theorem 3.1, it follows from
(3.15) that
$$
\begin{aligned}
&b\int_0^L\int_0^LG(x,s)w''(s)dsdx
\leq
2bE(a)D(a)\|p\|_1
\left[
1+c\int_0^L\int_0^LG(x,s)dsdx
\right]\\
&\quad\quad\quad\quad\quad\quad\quad\quad\quad\quad\quad\quad
<
1+c\int_0^L\int_0^LG(x,s)dsdx,
\end{aligned}
$$
where the last inequality follows from the equivalent form (3.3) of
hypothesis (H); hence (3.14) holds. Since the right-hand side of
(3.13) is nonnegative, it follows from (3.13) that
$$
\int_0^Lw(x)dx\geq0,
$$
and therefore \(w\in\mathcal X_+\). The uniqueness in \(\mathcal X_+\)
follows from Theorem 3.2. The proof is complete.\hfill{$\Box$}

\noindent{\bf Remark 3.2}\ For usual downward traffic loads $p\geq0$, Theorem 3.3 shows that the solution obtained in Theorem 3.1 is the unique solution in the physically meaningful class $\mathcal X_{+}$. In this case, the restriction on the solution can be formally relaxed to
$$
a+b\int_{0}^{L}w(x)dx>0.
$$
Indeed, let $w$ satisfy this condition and suppose $I:=\int_{0}^{L}w(x)dx<0$. Then $w''''-(a+bI)w''=p-cI>0$. Since $a+bI>0$, the positivity of the corresponding Green function yields $w>0$ in $(0,L)$, contradicting $I<0$. Hence every solution satisfying this formally relaxed condition still belongs to \(\mathcal X_+\) and therefore coincides with the solution obtained in Theorem 3.1. It should be noted, however, that global uniqueness in the strict mathematical sense cannot be asserted without a restriction of this kind, since the present analysis does not exclude additional negative-integral solutions of large magnitude whose associated fixed points lie outside $B[0,2\|p\|_1]$. For such possible solutions, the corresponding coefficient $a+b\int_{0}^{L} w(x)dx$ may reach or cross the first spectral value $\frac{-\pi^{2}}{L^{2}}$, leading to sign changes in the corresponding Green function and possible resonance. Their existence is a separate mathematical question, mainly of mathematical rather than engineering interest.

\noindent{\bf Remark 3.3} The key assumption in Theorem 3.1, Theorem 3.2 and Theorem 3.3 is the
explicit condition (H), which does not impose a separate smallness assumption on the main span length \(L\).  Therefore, the present results
improve the existence result and the uniqueness result for small solutions
in [38, Theorems 5.1 and 5.4] when \(h\) is taken as the classical
expression (1.3), where sufficiently small parameters \(L\) and \(p\) are
required. They also improve Wang's results [40, Theorem 4.2, Theorem 6.1
and Remark 6.2], where the conditions involve restrictions on \(L\) and
other structural parameters which, in practical suspension-bridge settings,
can be satisfied essentially only for short main spans and high deck
stiffness \(EI\). In contrast, the present results do not impose such
separate restrictions on the bridge parameters beyond the explicit
condition (H).

\noindent{\bf Remark 3.4}
For the original classical Melan problem (1.4), (1.2), the key condition (H) is equivalent to the following condition on the original bridge parameters
and the actual live load \(P\):
$$
\frac{
2H^3
\left[
1+\frac{E_cA_c}{L_c}q^2
\left(
\frac{L^3}{12H^3}
-\frac{EI\,L}{H^4}
+\frac{2(EI)^{3/2}}{H^{9/2}}
\tanh\!\left(\frac{L}{2}\sqrt{\frac{H}{EI}}\right)
\right)
\right]^2
}{
\frac{E_cA_c}{L_c}qL^2
\left[
1+\frac{E_cA_c}{L_c}q^2
\left(
\frac{5L^3}{24H^3}
-\frac{EI\,L}{H^4}
+\frac{(EI)^{3/2}}{H^{9/2}}
\left(
2-\frac{HL^2}{4EI}
\right)
\tanh\!\left(\frac{L}{2}\sqrt{\frac{H}{EI}}\right)
\right)
\right]
}
>
\|P\|_1 .
\eqno (3.16)
$$
The illustrative examples and discussion in Section 4 will show that, after customary lane-load reduction and vehicle-spacing effects are taken into account, the above condition (3.16) for (1.4),(1.2) is satisfied for practical suspension bridges under realistic small and moderate downward live loads $P$ representative of ordinary service traffic; hence the corresponding deflection curve is uniquely determined.

\section{Iterative method for classical Melan equation and its application to actual bridges}

%a fixed point iterative method is proposed to construct this unique solution with error estimates. The applicability of the proposed approach and the efficiency of the iterative algorithm are illustrated on two examples. The results show that the iterative method is highly accurate with fast convergence.

In this section we consider the iterative method for problem (3.1),(1.2),
namely the classical Melan equation (1.4) divided by \(EI\), with the
same boundary conditions (1.2). We also demonstrate the application of the
main results to two actual suspension bridges.

According to Theorem 3.1, under hypothesis (H), the operator \(A\) has a
unique fixed point \(\varphi\in B[0,2\|p\|_1]\), and hence problem
(3.1),(1.2) has the corresponding fixed-point solution. Moreover, when
\(p\geq0\), Theorem 3.3 shows that this solution is the unique solution
with nonnegative integral. To approximate the fixed-point solution, we may
choose any
$$
\varphi_0\in B[0,2\|p\|_1]
\eqno (4.1)
$$
as the initial function for iteration. By the definition of operator $A$, based on $\varphi_{m} (m = 0, 1, ...)$ we can obtain $w_{m}$ through (2.6); And consequently, we can get $\varphi_{m+1}$ by (3.4).
By the Banach contraction mapping principle, \(\varphi_m\) converges to
the unique fixed point \(\varphi\) of \(A\) at a geometric rate. Moreover,
the following error estimate holds:
$$
\|\varphi_m-\varphi\|_1
\leq
\frac{1}{1-q}\|\varphi_{m+1}-\varphi_m\|_1
\leq
\frac{q}{1-q}\|\varphi_m-\varphi_{m-1}\|_1
\leq\cdots
\leq
\frac{q^{m-1}}{1-q}\|\varphi_2-\varphi_1\|_1
\leq
\frac{q^m}{1-q}\|\varphi_1-\varphi_0\|_1,
\eqno (4.2)
$$
where \(q\) is defined in (3.9). Consequently, from (4.2) together with (3.10) we obtain the following convergence result of the iterative method:

\noindent{\bf Theorem 4.1}\
Suppose that the key condition (H) of Theorem 3.1 is satisfied. Then the iterative
sequence \(w_m\) generated by the operator \(A\) converges to the
fixed-point solution \(w\) of problem (3.1),(1.2) obtained in Theorem 3.1.
Furthermore, the error estimate is given by
$$
\aligned
\|w_m-w\|_1
&\leq K\|\varphi_m-\varphi\|_1
\leq
\frac{K}{1-q}\|\varphi_{m+1}-\varphi_m\|_1\\
&\leq
\frac{qK}{1-q}\|\varphi_m-\varphi_{m-1}\|_1
\leq\cdots
\leq
\frac{q^{m-1}K}{1-q}\|\varphi_2-\varphi_1\|_1
\leq
\frac{q^mK}{1-q}\|\varphi_1-\varphi_0\|_1.
\endaligned
\eqno (4.3)
$$
where
$$
K=
\frac{1+2c\left[\frac{L^3}{12a} - \frac{L}{a^2}
+\frac{2}{a^{5/2}} \tanh\left(\frac{\sqrt{a}L}{2}\right)\right]}
{1+c\left[\frac{L^3}{12a} - \frac{L}{a^2}
+\frac{2}{a^{5/2}} \tanh\left(\frac{\sqrt{a}L}{2}\right)\right]}
\frac{L^2}{8a}.
\eqno (4.4)
$$
If, in addition, \(p\geq0\), then this limit \(w\) is the unique solution
with nonnegative integral.

%\noindent{\bf Theorem 4.1} Suppose that all the conditions of Theorem 3.1 are satisfied, and let $\varphi_0$ be an arbitrary initial function in the closed ball $B[0,2\|p\|_1]$. For $m = 0,1,2,\dots$, define the sequence $\{\varphi_m\}$ by $\varphi_{m+1} = A\varphi_m$, where the operator $A$ is given by (3.2), and let $w_m$ be the function determined from $\varphi_m$ via the explicit formula (2.9). Then $\{\varphi_m\}$ converges to the unique fixed point $\varphi$ of $A$, and consequently $\{w_m\}$ converges to the unique solution $w$ of problem (3.1),(1.2). Moreover, the following error estimate holds:\[\|w_m - w\|_1 \le K\,\|\varphi_m - \varphi\|_1 \le K\,\frac{q^m}{1-q}\,\|\varphi_1 - \varphi_0\|_1,\] where $q$ is defined in (3.9) and $K$ is given in (3.10) by $$K=LG(\frac{L}{2},\frac{L}{2})\frac{1+2c\int^{L}_{0}\int^{L}_{0}G(x,s)dsdx}{1+c\int^{L}_{0}\int^{L}_{0}G(x,s)dsdx}.\eqno (4.4)$$

\noindent{\bf Example 4.1.} We consider the actual bridge discussed by Gazzola et al. [38], which has also been studied by Wollmann [13] and Shin et al. [18], taking the deformation effects of the bridge towers and side spans into account. Here, we focus exclusively on calculating the deflection of the main span stiffening girder using classical fundamental equation, without considering the deformation effects of the bridge towers and side spans. This serves as the basis for preliminary design calculations and for the subsequent determination of actual deflections via full-bridge modeling.  The main span $L$ of the suspension bridge is $460$m, and the sag-span ratio is $n=\frac{46}{460}=\frac{1}{10}$, combining these with other material parameters(see Wollmann [13], Gazzola et al. [38], and Shin et al. [18]) of the bridge, we transform Melan equation (1.4) into the form of (3.1) and can calculate the corresponding coefficients in (3.1) as follows:
$$
\aligned
&a=\frac{H}{EI}=\frac{97750}{5.7\times 10^{7}}=0.00171,\ \ b=\frac{q}{H}\frac{E_{c}A_{c}}{EI}\frac{1}{L_{c}}=0.001739\frac{3.6\times 10^{7}}{5.7\times 10^{7}}\frac{1}{472}=2.326\times 10^{-6},\\
&c=(\frac{q}{H})^2\frac{E_{c}A_{c}}{EI}\frac{1}{L_{c}}=\frac{q}{H}b=4.045\times10^{-9}.
\endaligned \eqno (4.5)
$$
From hypothesis (H), we can calculate the following:
$$
\frac{2a^2
\left[1+c\left(\frac{L^3}{12a}-\frac{L}{a^2}+\frac{2}{a^{5/2}}\tanh\!\left(\frac{\sqrt aL}{2}\right)\right)\right]^2}{
bL^2\left[
1+
c\left(
\frac{5L^3}{24a}
-\frac{L}{a^2}
+\frac{1}{a^{5/2}}
\left(
2-\frac{aL^2}{4}
\right)
\tanh\!\left(\frac{\sqrt aL}{2}\right)
\right)
\right]}=1.00787\times10^{-4},
\eqno (4.6)$$
that is, the function $p$ in (3.1) should satisfy $\|p\|_{1}<1.00787\times10^{-4}$
and then the real live load $P$ on the bridge in (1.4) should meet the following condition:
$$
\|P\|_{1}<1.00787\times10^{-4}\times EI=1.00787\times10^{-4}\times5.7\times 10^{7}=5744.859(\mathrm{kN}).
\eqno (4.7)
$$Now, without loss of generality, we consider the same two live load cases as in [38]: a coach of length $10$m with a weight density of $10$ kN/m, and a freight train of length $230m$ with a weight density of $20$ kN/m, both centered at midspan.

\noindent{\bf Case1.} A coach of length $10$m with a weight density of $10$ kN/m in the middle of the span. The actual live load function $P(x)$ in (1.4) for this case is:
$$P(x)=
\begin{cases}
0, & 0\leq x\leq 225,\\
10, & 225< x\leq 235,\\
0, & 235 < x\leq 460.
\end{cases}\eqno (4.8)$$
Since the total actual live load is $\|P\|_1=\int_0^{460} |P(x)|dx=100(\mathrm{kN})$, it follows from (4.7) that the hypothesis (H) of Theorem 3.1 holds in this case.
By $EI=5.7\times 10^{7}$, the function $p(x)$ in equation (3.1) is
$$p(x)=\begin{cases}
0, & 0\leq x\leq 225,\\
\frac{10}{5.7\times 10^{7}}, & 225< x\leq 235,\\
0, & 235 < x\leq 460.
\end{cases}\eqno (4.9)$$
For this normalized load and the parameters \(a,b,c\) given in (4.5), the problem corresponding to (3.1) becomes
$$
\aligned
&w''''(x)-(a+b\int_0^{460} w(x)dx)w''(x)+c\int_0^{460} w(x)dx=p(x),\ \ \ \ x\in (0,460),\\
&w(0)=w(460)=w''(0)=w''(460)=0.\\
\endaligned \eqno (4.10)
$$

According to the iterative method and (4.1), for computational convenience, we may choose
$$
\varphi_{0}\equiv0
$$
as the initial function for iteration. Substituting $\varphi_{0}\equiv0$ into the right side of (2.6) we can obtain
$$
w_{0}\equiv0.
$$
Then by (3.4) we have
$$
\varphi_{1}(x)=p(x).
$$
Using the explicit form of $p(x)$ in (4.9), we obtain the following expression for $w_{1}$ after substitution and calculation:
$$
w_{1}(x)=p[\int^{235 }_{225}G(x,s)ds-\frac{c\int^{460}_{0}G(x,s)ds}{1+c\int^{460}_{0}\int^{460}_{0}G(x,s)dsdx}\int^{460}_{0}\int^{235}_{225}G(x,s)dsdx],
\eqno (4.11)$$
where $p=\frac{10}{5.7\times 10^{7}}$.
Since by (2.5) we can calculate the following needed integrals:
$$
\aligned
\int^{235}_{225}G(x,s)\,ds=
\left\{
\begin{aligned}
&
\frac{5}{a}x
+\frac{\sinh(\sqrt a x)
\left(
\cosh(225\sqrt a)-\cosh(235\sqrt a)
\right)}
{a^{2}\sinh(460\sqrt a)},
\quad \ \ \ \ \ \ \ \ \ \ \ \ \ \ \ \ \ \  0\leq x\leq225,
\\[2mm]
&
\frac{460x-x^{2}-225^{2}}{2a}
-\frac{1}{a^{2}}
+\frac{
\cosh(225\sqrt a)
\left(
\sinh(\sqrt a(460-x))+\sinh(\sqrt a x)
\right)}
{a^{2}\sinh(460\sqrt a)},
\\
&\hspace{10.8cm}225<x\leq235,
\\[2mm]
&
\frac{5}{a}(460-x)
+\frac{\sinh(\sqrt a(460-x))
\left(
\cosh(225\sqrt a)-\cosh(235\sqrt a)
\right)}
{a^{2}\sinh(460\sqrt a)},
\  235<x\leq460,
\end{aligned}
\right.
\endaligned
\eqno (4.12)
$$
$$
\aligned
\int^{460}_{0}G(x,s)ds=& \frac{460x-x^2}{2a}+\frac{\sinh(\sqrt{a}x)+\sinh(\sqrt{a}(460-x))-\sinh(\sqrt{a}460)}{a^2\sinh(\sqrt{a}460)},\\
\endaligned \eqno (4.13)
$$
$$
\int^{460}_{0}\int^{235}_{225}G(x,s)dsdx\approx1.512347\times10^{8},
$$
$$\int^{460}_{0}\int^{460}_{0}G(x,s)dsdx\approx4.602697\times10^{9},\eqno (4.14)$$
then by (4.11) we have
\begin{small}
$$
\aligned
&w_{1}(x)
=p\left[
\int^{235}_{225}G(x,s)\,ds
-k\int^{460}_{0}G(x,s)\,ds
\right]
\\
&=p\left\{
\begin{array}{@{}l@{}}
\displaystyle
\frac{x(10-460k)+kx^{2}}{2a}
+\frac{1}{a^{2}\sinh(460\sqrt a)}
\bigg[
\sinh(\sqrt a x)
\left(\cosh(225\sqrt a)-\cosh(235\sqrt a)-k\right)
\\
\displaystyle
\hspace{2.0em}
-k\sinh(\sqrt a(460-x))
+k\sinh(460\sqrt a)
\bigg],
\hspace{3.2em}0\leq x\leq225,
\\[4mm]
\displaystyle
\frac{460x(1-k)+(k-1)x^{2}-225^{2}}{2a}
+\frac{1}{a^{2}\sinh(460\sqrt a)}
\bigg[
\left(\sinh(\sqrt a(460-x))+\sinh(\sqrt a x)\right)
\left(\cosh(225\sqrt a)-k\right)
\\
\displaystyle
\hspace{2.0em}
-(1-k)\sinh(460\sqrt a)
\bigg],
\hspace{9.9em}225<x\leq235,
\\[4mm]
\displaystyle
\frac{4600-(10+460k)x+kx^{2}}{2a}
+\frac{1}{a^{2}\sinh(460\sqrt a)}
\bigg[
\sinh(\sqrt a(460-x))
\left(\cosh(225\sqrt a)-\cosh(235\sqrt a)-k\right)
\\
\displaystyle
\hspace{2.0em}
-k\sinh(\sqrt a x)
+k\sinh(460\sqrt a)
\bigg],
\hspace{6.3em}235<x\leq460.
\end{array}
\right.
\endaligned
\eqno (4.15)
$$
\end{small}
where $p=\frac{10}{5.7\times 10^{7}}$, $k=\frac{c\int^{460}_{0}\int^{235}_{225}G(x,s)dsdx}{1+c\int^{460}_{0}\int^{460}_{0}G(x,s)dsdx}=0.03118295, a=0.00171.$ It should be noted that although the function $w_{1}$ appears piecewise, it can be verified that it is twice continuously differentiable and symmetric with respect to the center point $x=230$.

Then by (3.4) and calculation we have
\begin{small} $$ \aligned &\varphi_{2}(x) =p(x)+b\int_{0}^{460}w_{1}(x)\,dx\,w_{1}''(x)=p(x)+2.326\times10^{-6}\times1.352458\,w_{1}''(x) \\ &=
\left\{
\begin{array}{@{}l@{}}
\displaystyle
3.145817\times10^{-6}\frac{p}{a}
\left[
k+
\frac{
\left(\cosh(225\sqrt a)-\cosh(235\sqrt a)-k\right)\sinh(\sqrt a x)
-k\sinh(\sqrt a(460-x))}
{\sinh(460\sqrt a)}
\right],
\\[-1mm]
\hfill\hspace{2em}0\leq x\leq225,
\\[3mm]
\displaystyle
\frac{10}{5.7\times10^{7}}
+
3.145817\times10^{-6}\frac{p}{a}
\left[
k-1+
\frac{
\left(\cosh(225\sqrt a)-k\right)
\left(\sinh(\sqrt a x)+\sinh(\sqrt a(460-x))\right)}
{\sinh(460\sqrt a)}
\right],
\\[-1mm]
\hfill\hspace{2em}225<x\leq235,
\\[3mm]
\displaystyle
3.145817\times10^{-6}\frac{p}{a}
\left[
k+
\frac{
\left(\cosh(225\sqrt a)-\cosh(235\sqrt a)-k\right)\sinh(\sqrt a(460-x))
-k\sinh(\sqrt a x)}
{\sinh(460\sqrt a)}
\right],
\\[-1mm]
\hfill\hspace{2em}235<x\leq460.
\end{array}
\right.
\\[2mm]&\approx
\left\{
\begin{array}{@{}l@{\hspace{2.5em}}r@{}}
\displaystyle
1.006421\times10^{-11}
-9.951379\times10^{-15}\sinh(\sqrt{0.00171}\,x)
&
\\[-1mm]
\displaystyle
\ \ \ \ \ \ -1.103229\times10^{-19}\sinh(\sqrt{0.00171}(460-x)),&0\leq x\leq225,
\\[3mm]
\displaystyle
1.751259\times10^{-7}
+1.94309\times10^{-14}\sinh(\sqrt{0.00171}\,x)
&
\\[-1mm]
\displaystyle
\ \ \ \ \ \ +1.94309\times10^{-14}\sinh(\sqrt{0.00171}(460-x)),&225<x\leq235,
\\[3mm]
\displaystyle
1.006421\times10^{-11}
-1.103229\times10^{-19}\sinh(\sqrt{0.00171}\,x)
&
\\[-1mm]
\displaystyle
\ \ \ \ \ \ -9.951379\times10^{-15}\sinh(\sqrt{0.00171}(460-x)),&235<x\leq460.
\end{array}
\right. \endaligned \eqno(4.16) $$ \end{small}
To obtain $w_{2}$ by substituting $\varphi_{2}$ into the right-hand side of (2.6), we first compute the following:
\begin{small} $$ \aligned &\int^{460}_{0}G(x,s)\varphi_{2}(s)\,ds \\ &\approx \left\{ \begin{aligned} & -6.883629\times10^{-6} +5.133883\times10^{-4}x -2.942751\times10^{-9}x^{2} \\ &\quad -1.828638\times10^{-6}\sinh(\sqrt{0.00171}x) +7.545766\times10^{-14}\sinh(\sqrt{0.00171}(460-x)) \\ &\quad -7.036539\times10^{-11}x\sinh(\sqrt{0.00171}x) +7.793127\times10^{-16}x\sinh(\sqrt{0.00171}(460-x)),0\leq x\leq225, \\ & -2.651959 +2.355495\times10^{-2}x -5.120641\times10^{-5}x^{2} \\ &\quad +3.566573\times10^{-6}\sinh(\sqrt{0.00171}x) +3.629774\times10^{-6}\sinh(\sqrt{0.00171}(460-x)) \\ &\quad +1.373943\times10^{-10}x\sinh(\sqrt{0.00171}x) -1.373943\times10^{-10}x\sinh(\sqrt{0.00171}(460-x)),225<x\leq235,\\ & 2.355291\times10^{-1} -5.10681\times10^{-4}x -2.942751\times10^{-9}x^{2} \\ &\quad +4.339415\times10^{-13}\sinh(\sqrt{0.00171}x) -1.861006\times10^{-6}\sinh(\sqrt{0.00171}(460-x)) \\ &\quad -7.793127\times10^{-16}x\sinh(\sqrt{0.00171}x) +7.036539\times10^{-11}x\sinh(\sqrt{0.00171}(460-x)),235<x\leq460.
\end{aligned} \right.
 \endaligned \eqno (4.17) $$ \end{small}
and then we can calculate the following result:
$$\int^{460}_{0}\int^{460}_{0}G(x,s)\varphi_{2}(s)dsdx\approx 26.52968.\eqno (4.18)$$
Now, substituting (4.13),(4.14),(4.17) and (4.18) into (2.6), we can obtain $w_{2}$ as follows:
\begin{small}
$$
\aligned
&w_{2}(x)=\int^{L}_{0}G(x,s)\varphi_{2}(s)ds-
\frac{c\int^{L}_{0}\int^{L}_{0}G(x,s)\varphi_{2}(s)dsdx}{1+c\int^{L}_{0}\int^{L}_{0}G(x,s)dsdx}\int^{L}_{0}G(x,s)ds\\
&=\int^{L}_{0}G(x,s)\varphi_{2}(s)ds-5.470132\times 10^{-9}
\int^{L}_{0}G(x,s)ds\\&\approx\left\{ \begin{aligned} & 1.863822\times10^{-3} -2.223604\times10^{-4}x +1.596511\times10^{-6}x^{2} \\ &\quad -1.828658\times10^{-6}\sinh(\sqrt{0.00171}x) -2.043104\times10^{-11}\sinh(\sqrt{0.00171}(460-x)) \\ &\quad -7.036539\times10^{-11}x\sinh(\sqrt{0.00171}x) +7.793127\times10^{-16}x\sinh(\sqrt{0.00171}(460-x)),\ 0\leq x\leq225, \\[2mm] & -2.650088 +2.28192\times10^{-2}x -4.960695\times10^{-5}x^{2} \\ &\quad +3.566552\times10^{-6}\sinh(\sqrt{0.00171}x) +3.629754\times10^{-6}\sinh(\sqrt{0.00171}(460-x)) \\ &\quad +1.373943\times10^{-10}x\sinh(\sqrt{0.00171}x) -1.373943\times10^{-10}x\sinh(\sqrt{0.00171}(460-x)),\ 225<x\leq235, \\[2mm] & 2.373998\times10^{-1} -1.24643\times10^{-3}x +1.596511\times10^{-6}x^{2} \\ &\quad -2.007255\times10^{-11}\sinh(\sqrt{0.00171}x) -1.861026\times10^{-6}\sinh(\sqrt{0.00171}(460-x)) \\ &\quad -7.793127\times10^{-16}x\sinh(\sqrt{0.00171}x) +7.036539\times10^{-11}x\sinh(\sqrt{0.00171}(460-x)),\ 235<x\leq460. \end{aligned} \right.
\endaligned\eqno (4.19)
$$
\end{small}
Then by (3.4) and calculation we have
{\small\begin{equation}
\begin{aligned}
&\varphi_{3}(x)
=p(x)+b\int_{0}^{460}w_{2}(x)\,dx\,w_{2}''(x)=p(x)+2.326\times10^{-6}\times1.352319\,w_{2}''(x)
\\
&\approx
\left\{
\begin{aligned}&
1.004363\times10^{-11}
-9.83598\times10^{-15}\sinh(\sqrt{0.00171}x)
-1.098944\times10^{-19}\sinh(\sqrt{0.00171}(460-x))
\\&\quad
-3.784811\times10^{-19}x\sinh(\sqrt{0.00171}x)
+4.191764\times10^{-24}x\sinh(\sqrt{0.00171}(460-x))
\\
&\quad
-1.830527\times10^{-17}\cosh(\sqrt{0.00171}x)
-2.02735\times10^{-22}\cosh(\sqrt{0.00171}(460-x)),
\quad 0\leq x\leq225,
\\[2mm]&
1.751265\times10^{-7}
+1.918376\times10^{-14}\sinh(\sqrt{0.00171}x)
+1.95237\times10^{-14}\sinh(\sqrt{0.00171}(460-x))
\\&\quad
+7.39016\times10^{-19}x\sinh(\sqrt{0.00171}x)
-7.39016\times10^{-19}x\sinh(\sqrt{0.00171}(460-x))
\\&\quad
+3.574257\times10^{-17}\cosh(\sqrt{0.00171}x)
+3.574257\times10^{-17}\cosh(\sqrt{0.00171}(460-x)),
\quad 225<x\leq235,
\\[2mm]&
1.004363\times10^{-11}
-1.079662\times10^{-19}\sinh(\sqrt{0.00171}x)
-1.001008\times10^{-14}\sinh(\sqrt{0.00171}(460-x))
\\&\quad
-4.191764\times10^{-24}x\sinh(\sqrt{0.00171}x)
+3.784811\times10^{-19}x\sinh(\sqrt{0.00171}(460-x))
\\&\quad
-2.02735\times10^{-22}\cosh(\sqrt{0.00171}x)
-1.830527\times10^{-17}\cosh(\sqrt{0.00171}(460-x)),
\quad 235<x\leq460.
\end{aligned}\right.
\end{aligned}\tag{4.20}\end{equation}}
Because of the tedious computation involved in evaluating $\int^{L}_{0}G(x,s)\varphi_{3}(s)ds$, we omit the derivation of the expression for $w_{3}$. The results of the two iterations are shown in Figure 2 below drawn with MATLAB:
\begin{figure}[H]
\centering
\includegraphics[width=0.9\linewidth]{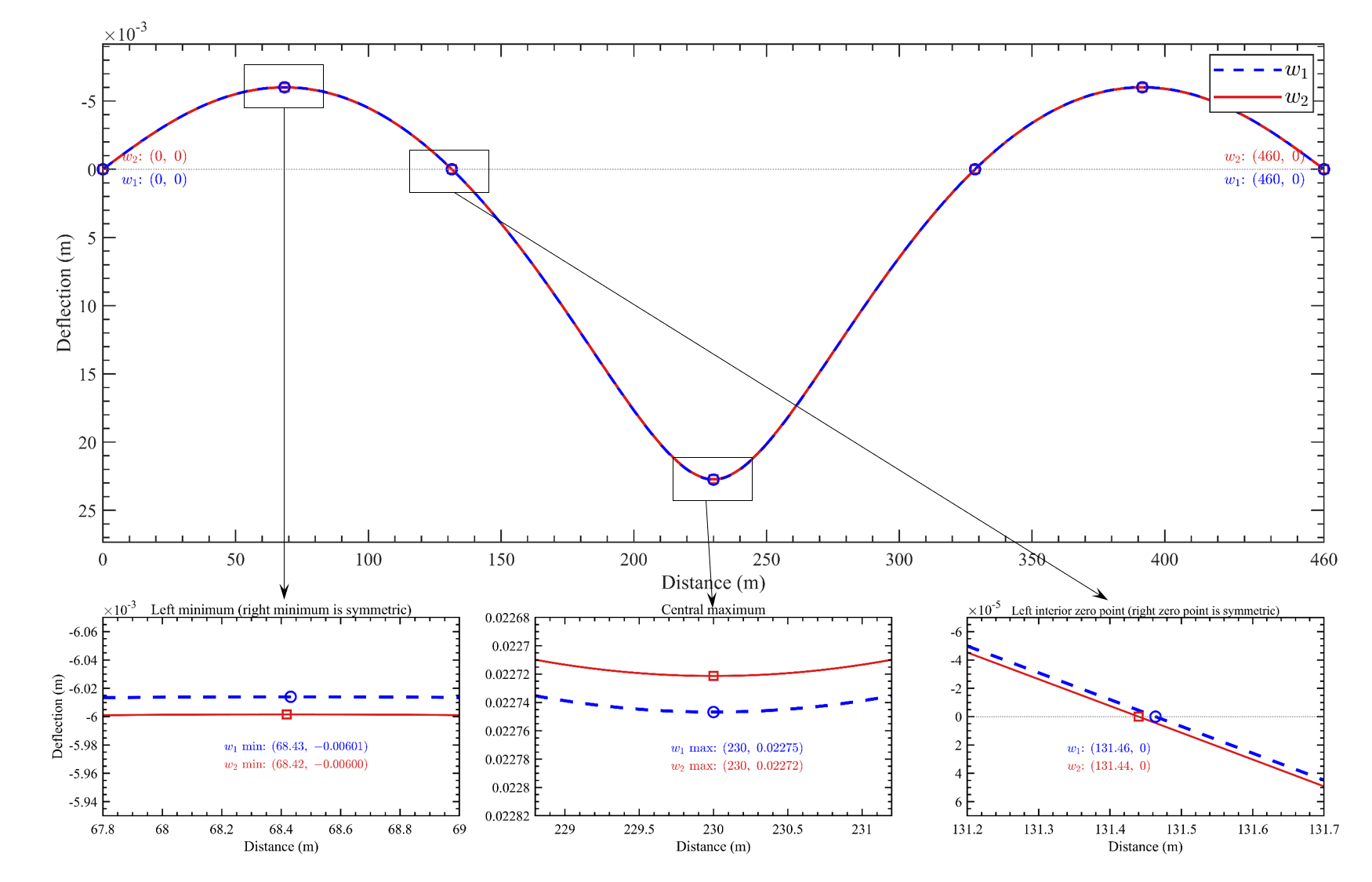}
\caption{The first two iterates \(w_1\) and \(w_2\) for the suspension
bridge with \(L=460\,{\rm m}\) under a uniformly distributed live load of
\(10\,{\rm kN/m}\) on \(225\le x\le235\).}
\end{figure}
As will be quantified by the error
estimates below, even the first iterate $w_{1}$ already gives a good approximation to the exact solution $w$ of (4.10), while the second
iterate $w_{2}$ is remarkably close to it. We next quantify this by
investigating the errors of the first two iterates using the theoretical global contraction constant.

According to the load function (4.9) for this case, we have
$$
\|p\|_{1}=\|\varphi_{1}-\varphi_{0}\|_{1}=\int_0^{460}|p(x)|dx=1.754386\times10^{-6}.
\eqno (4.21)
$$
then, using (3.9) and (4.6) we can calculate that the corresponding global contraction constant is
$$q=\|p\|_{1}\times\frac{1}{1.00787\times10^{-4}}=1.754386\times10^{-6}\times\frac{1}{1.00787\times10^{-4}}=0.01740687.$$
With this and the iterates $\varphi_{1},\varphi_{2},\varphi_{3}$, (4.2) yields the following a posteriori theoretical error estimate between the iterative fixed points and the true fixed point:
$$
\|\varphi_1-\varphi\|_1\leq\frac{1}{1-q}\|\varphi_2-\varphi_1\|_1=
\frac{1}{1-0.01740687}\times 4.549183\times10^{-9}=4.629773\times10^{-9},$$
$$
\aligned
\|\varphi_2-\varphi\|_1&\leq \frac{1}{1-q}\|\varphi_{3}-\varphi_2\|_1=\frac{1}{1-0.01740687}\times7.088715\times10^{-12}=7.214293\times10^{-12}.\\\endaligned$$Consequently, combining the above results with the value of $K=3.014\times10^{7}$ computed from (4.4), Theorem 4.1 yields the following estimate for the iterative solutions of the problem (4.10):
$$\begin{aligned}\|w_1-w\|_1\leq K\|\varphi_1-\varphi\|_1\leq3.014\times10^{7}\times4.629773\times10^{-9}=
0.1395414,
\end{aligned}
$$
$$
\aligned
\|w_{2}-w\|_{1}\leq K\|\varphi_2-\varphi\|_1\leq3.014\times10^{7}\times7.214293\times10^{-12}=0.0002174388.\endaligned$$
By the above formulas, the mean absolute errors of $w_{1}$ and $w_{2}$(relative to the exact solution $w$) over the entire main span do not exceed $0.0003033508$ and
$4.72693\times10^{-7} m,$ respectively. Comparing the above results with the amplitudes of $w_{1}$ and $w_{2}$ computed from (4.15) and (4.19):
$$\underset{x\in (0,460)}\max w_{1}(x)=w_{1}(230)=0.02275 m,\ \underset{x\in (0,460)}\min w_{1}(x)=w_{1}(68.43185)=w_{1}(391.56815)=-0.00601 m, $$
$$\underset{x\in (0,460)}\max w_{2}(x)=w_{2}(230)=0.02272 m,\ \underset{x\in (0,460)}\min w_{2}(x)=w_{2}(68.41828)=w_{2}(391.58172)=-0.006 m $$
we see that the mean absolute error of \(w_1\) is already much smaller than the displacement amplitude, and $w_{2}$ is even more accurate. Although the above global error estimates are already sufficiently sharp, the global contraction constant still gives a conservative description of the actual convergence behaviour, as indicated by the following observed ratios.  Indeed, by direct calculation we obtain
$$
\frac{\|\varphi_2-\varphi_1\|_1} {\|\varphi_1-\varphi_0\|_1} \approx \frac{4.549183\times10^{-9}}{1.754386\times10^{-6}}\approx2.59304\times10^{-3}:=q_{1}, $$
$$ \frac{\|\varphi_3-\varphi_2\|_1} {\|\varphi_2-\varphi_1\|_1} \approx \frac{7.088715\times10^{-12}}{4.549183\times10^{-9}}\approx 1.55824\times10^{-3}:=q_{2}.
$$The above ratios suggest that the convergence along the computed iteration
is much faster than that indicated by the global contraction constant; more
precisely, the observed ratio is nearly one order of magnitude
smaller than the global contraction constant. This motivates us to
introduce the observed local contraction factor
$$
q_{\rm loc}:=\max\{q_1,q_2\}
=2.59304\times10^{-3}
<q=0.01740687,
$$
which, when used as the effective contraction factor, would lead to faster
and smaller a posteriori error estimates that better reflect the actual
numerical convergence of the iteration.  Since the global estimates are
already sufficiently small here, we do not list these local estimates in
detail.

\noindent{\bf Case2.} A freight train of length $230m$ with a weight density of $20$ kN/m in the middle of the span. The actual live load function $P(x)$ in (1.4) for this case is:
$$P(x)=
\begin{cases}
0, & 0\leq x\leq 115,\\
20, & 115< x\leq 345,\\
0, & 345 < x\leq 460.
\end{cases}$$
and the corresponding function $p$ in equation (4.10) is
$$p(x)=
\begin{cases}
0, & 0\leq x\leq 115,\\
\frac{20}{5.7\times 10^{7}}, & 115< x\leq 345,\\
0, & 345 < x\leq 460.
\end{cases}\eqno (4.22)$$
Since $\|P\|_{1}=\int_0^{460} |P(x)|dx=4600(\mathrm{kN})$, it follows from (4.7) that the Condition (H) of Theorem 3.1 holds under this live load.
Following the same procedure as in Case 1, we take $\varphi_{0}\equiv0$
as the initial function and then, by (2.6) and (3.4), successively compute the following:
$$
\aligned
&w_{0}\equiv0;\\
&\varphi_{1}(x)=p(x);\\
\endaligned
$$
\begin{small}
$$
\aligned
&w_{1}(x)=p[\int^{345 }_{115}G(x,s)ds-\frac{c\int^{460}_{0}G(x,s)ds}{1+c\int^{460}_{0}\int^{460}_{0}G(x,s)dsdx}\int^{460}_{0}\int^{345}_{115}G(x,s)dsdx]\\
&=p\left\{
\begin{array}{@{}l@{}}
\displaystyle
\frac{460x(1-2k)+2kx^{2}}{4a}
+\frac{1}{a^{2}\sinh(460\sqrt a)}
\bigg[
\sinh(\sqrt a x)
\left(
\cosh(115\sqrt a)-\cosh(345\sqrt a)-k
\right)
\\
\displaystyle
\hspace{2.0em}
-k\sinh(\sqrt a(460-x))
+k\sinh(460\sqrt a)
\bigg],
\hspace{3.2em}0\leq x\leq115,
\\[4mm]
\displaystyle
\frac{460x(1-k)+(k-1)x^{2}-\frac{460^{2}}{16}}{2a}
+\frac{1}{a^{2}\sinh(460\sqrt a)}
\bigg[
\left(\sinh(\sqrt a(460-x))+\sinh(\sqrt a x)\right)
\left(\cosh(115\sqrt a)-k\right)
\\
\displaystyle
\hspace{2.0em}
-(1-k)\sinh(460\sqrt a)
\bigg],
\hspace{9.3em}115<x\leq345,
\\[4mm]
\displaystyle
\frac{460^{2}-460x(1+2k)+2kx^{2}}{4a}
+\frac{1}{a^{2}\sinh(460\sqrt a)}
\bigg[
\sinh(\sqrt a(460-x))
\left(
\cosh(115\sqrt a)-\cosh(345\sqrt a)-k
\right)
\\
\displaystyle
\hspace{2.0em}
-k\sinh(\sqrt a x)
+k\sinh(460\sqrt a)
\bigg],
\hspace{6.3em}345<x\leq460,
\end{array}
\right.
\endaligned\eqno (4.23)
$$
\end{small}
\nopagebreak[4]
\noindent where $p=\frac{20}{5.7\times 10^{7}}$, $k=\frac{c\int^{460}_{0}\int^{345}_{115}G(x,s)dsdx}{1+c\int^{460}_{0}\int^{460}_{0}G(x,s)dsdx}=0.656222$, and $a=0.00171$;
\par\pagebreak[3]
\begin{small}
$$
\aligned
&\varphi_{2}(x)
=p(x)+b\int_{0}^{460}w_{1}(x)\,dx\,w_{1}''(x)
=p(x)+2.326\times10^{-6}\times56.92295\,w_{1}''(x)
\\
&=
\left\{
\begin{array}{@{}l@{}}
\displaystyle
0.0001324028\frac{p}{a}
\left[
k+
\frac{
\left(
\cosh\left(\frac{460\sqrt a}{4}\right)
-\cosh\left(\frac{3\cdot460\sqrt a}{4}\right)
-k
\right)\sinh(\sqrt a x)
-k\sinh(\sqrt a(460-x))}
{\sinh(460\sqrt a)}
\right],
\\[-1mm]
\hfill\hspace{2em}0\leq x\leq115,
\\[3mm]
\displaystyle
\frac{20}{5.7\times10^{7}}
+
0.0001324028\frac{p}{a}
\left[
k-1+
\frac{
\left(
\cosh\left(\frac{460\sqrt a}{4}\right)-k
\right)
\left(
\sinh(\sqrt a x)+\sinh(\sqrt a(460-x))
\right)}
{\sinh(460\sqrt a)}
\right],
\\[-1mm]
\hfill\hspace{2em}115<x\leq345,
\\[3mm]
\displaystyle
0.0001324028\frac{p}{a}
\left[
k+
\frac{
\left(
\cosh\left(\frac{460\sqrt a}{4}\right)
-\cosh\left(\frac{3\cdot460\sqrt a}{4}\right)
-k
\right)\sinh(\sqrt a(460-x))
-k\sinh(\sqrt a x)}
{\sinh(460\sqrt a)}
\right],
\\[-1mm]
\hfill\hspace{2em}345<x\leq460
\end{array}
\right.
\\[2mm]&\approx
\left\{
\begin{array}{@{}l@{\hspace{2.5em}}r@{}}
\displaystyle
1.782818\times10^{-8}
-2.337428\times10^{-10}\sinh(\sqrt{0.00171}\,x)
&
\\[-1mm]
\displaystyle
\ \ \ \ \
-1.954307\times10^{-16}\sinh(\sqrt{0.00171}(460-x)),
&0\leq x\leq115,
\\[3mm]
\displaystyle
3.415375\times10^{-7}
+1.71119\times10^{-14}\sinh(\sqrt{0.00171}\,x)
&
\\[-1mm]
\displaystyle
\ \ \ \ \
+1.71119\times10^{-14}\sinh(\sqrt{0.00171}(460-x)),
&115<x\leq345,
\\[3mm]
\displaystyle
1.782818\times10^{-8}
-1.954307\times10^{-16}\sinh(\sqrt{0.00171}\,x)
&
\\[-1mm]
\displaystyle
\ \ \ \ \
-2.337428\times10^{-10}\sinh(\sqrt{0.00171}(460-x)),
&345<x\leq460;
\end{array}
\right.
\endaligned
\eqno (4.24)
$$
\end{small}
\begin{small}
$$
\aligned
&w_{2}(x)=\int^{L}_{0}G(x,s)\varphi_{2}(s)ds-
\frac{c\int^{L}_{0}\int^{L}_{0}G(x,s)\varphi_{2}(s)dsdx}{1+c\int^{L}_{0}\int^{L}_{0}G(x,s)dsdx}\int^{L}_{0}G(x,s)ds\\
&=\int^{L}_{0}G(x,s)\varphi_{2}(s)ds-2.29295\times10^{-7}
\int^{L}_{0}G(x,s)ds\\
&\approx\left\{
\begin{aligned}
&
0.06622162
-0.006672984x
+0.00006183241x^{2}
\\
&\quad
-0.0006824781\sinh(\sqrt{0.00171}x)
-7.25958\times10^{-10}\sinh(\sqrt{0.00171}(460-x))
\\
&\quad
-1.652777\times10^{-6}x\sinh(\sqrt{0.00171}x)
+1.363757\times10^{-12}x\sinh(\sqrt{0.00171}(460-x)),
 0\leq x\leq115,
\\[2mm]
&
-1.286962
+0.01509694x
-0.00003281943x^{2}
\\
&\quad
+2.22921\times10^{-8}\sinh(\sqrt{0.00171}x)
+7.795067\times10^{-8}\sinh(\sqrt{0.00171}(460-x))
\\
&\quad
+1.209969\times10^{-10}x\sinh(\sqrt{0.00171}x)
-1.209969\times10^{-10}x\sinh(\sqrt{0.00171}(460-x)),
115<x\leq345,
\\[2mm]
&
10.08039
-0.05021283x
+0.00006183241x^{2}
\\
&\quad
-9.862977\times10^{-11}\sinh(\sqrt{0.00171}x)
-0.001442755\sinh(\sqrt{0.00171}(460-x))
\\
&\quad
-1.363757\times10^{-12}x\sinh(\sqrt{0.00171}x)
+1.652777\times10^{-6}x\sinh(\sqrt{0.00171}(460-x)),
345<x\leq460;
\end{aligned}
\right.
\endaligned\eqno (4.25)
$$
\end{small}
{\small \begin{equation} \begin{aligned} &\varphi_{3}(x)
=p(x)+b\int_{0}^{460}w_{2}(x)\,dx\,w_{2}''(x)=p(x)+2.326\times10^{-6}\times56.68584\,w_{2}''(x) \\ &\approx \left\{ \begin{aligned} & 1.630536\times10^{-8} -1.538754\times10^{-10}\sinh(\sqrt{0.00171}x) -1.636786\times10^{-16}\sinh(\sqrt{0.00171}(460-x)) \\ &\quad -3.726443\times10^{-13}x\sinh(\sqrt{0.00171}x) +3.074804\times10^{-19}x\sinh(\sqrt{0.00171}(460-x)) \\ &\quad -1.802297\times10^{-11}\cosh(\sqrt{0.00171}x) -1.487131\times10^{-17}\cosh(\sqrt{0.00171}(460-x)), \quad 0\leq x\leq115, \\[2mm] & 3.422226\times10^{-7} +5.026103\times10^{-15}\sinh(\sqrt{0.00171}x) +1.75752\times10^{-14}\sinh(\sqrt{0.00171}(460-x)) \\ &\quad +2.728064\times10^{-17}x\sinh(\sqrt{0.00171}x) -2.728064\times10^{-17}x\sinh(\sqrt{0.00171}(460-x)) \\ &\quad +1.31943\times10^{-15}\cosh(\sqrt{0.00171}x) +1.31943\times10^{-15}\cosh(\sqrt{0.00171}(460-x)), \quad 115<x\leq345, \\[2mm] & 1.630536\times10^{-8} -2.223763\times10^{-17}\sinh(\sqrt{0.00171}x) -3.252918\times10^{-10}\sinh(\sqrt{0.00171}(460-x)) \\ &\quad -3.074804\times10^{-19}x\sinh(\sqrt{0.00171}x) +3.726443\times10^{-13}x\sinh(\sqrt{0.00171}(460-x)) \\ &\quad -1.487131\times10^{-17}\cosh(\sqrt{0.00171}x) -1.802297\times10^{-11}\cosh(\sqrt{0.00171}(460-x)), \quad 345<x\leq460. \end{aligned} \right. \end{aligned} \tag{4.26} \end{equation} }
Similarly to Case 1, we omit the explicit expression for
$w_{3}$ owing to the tedious integration involved. The results of the two iterations are shown in Figure 3 below drawn with MATLAB:
\begin{figure}[H]
\centering
\includegraphics[width=0.9\linewidth]{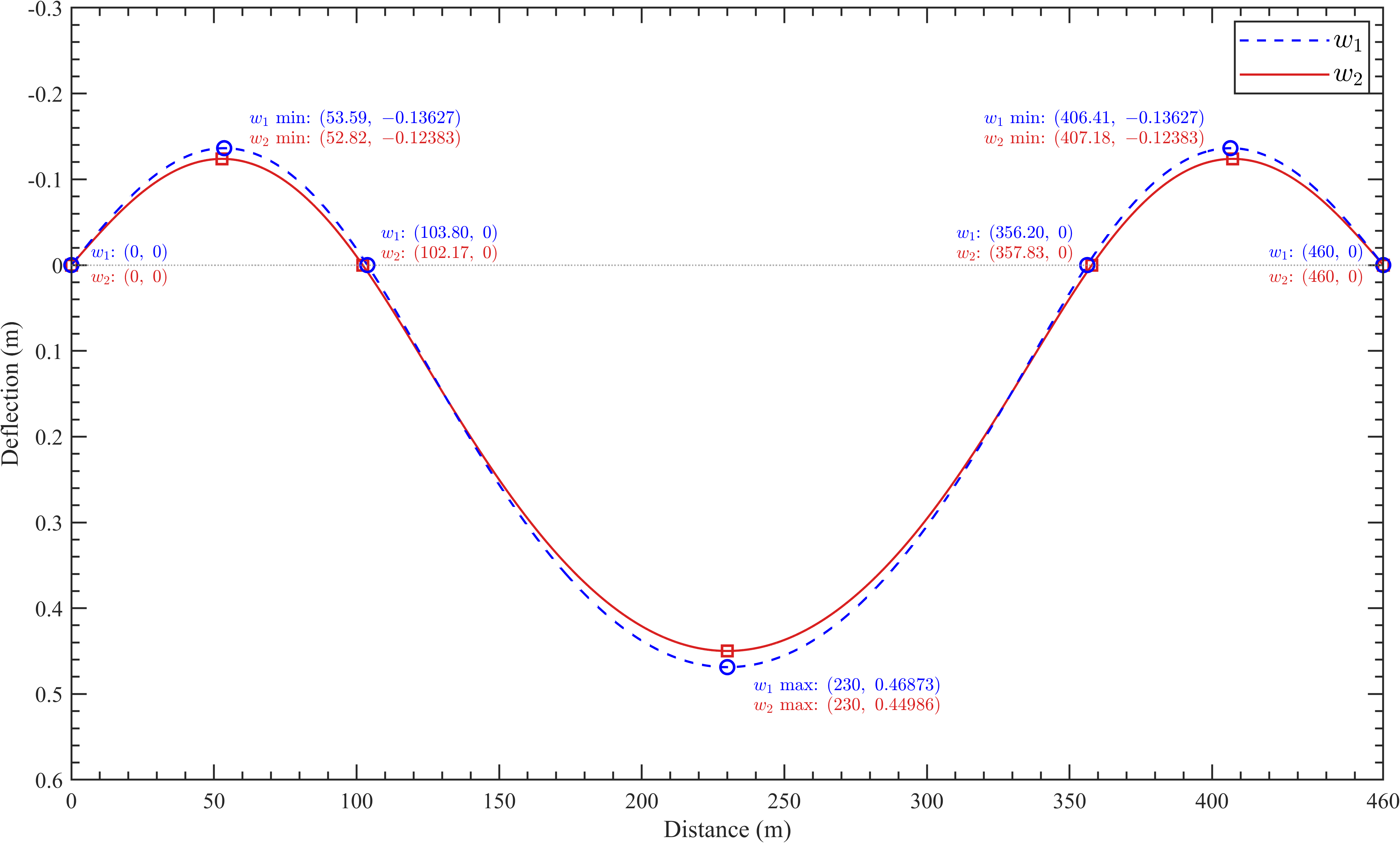}
\caption{The first two iterates \(w_1\) and \(w_2\) for the suspension
bridge with \(L=460\,{\rm m}\) under a uniformly distributed live load of
\(20\,{\rm kN/m}\) on \(115\le x\le345\).}
\end{figure}
As will be quantified by the error
estimates below, although the global contraction estimate gives
relatively conservative error bounds in this case, the estimate based on
the observed local contraction factor shows that $w_{1}$ is still an
acceptable approximation, while $w_{2}$ provides a considerably more
accurate approximation. We next quantify this by investigating the errors of the first two iterates with respect to the exact solution, using both the theoretical global contraction constant and the observed local contraction factor. According to the load function (4.22), we have
$$
\|p\|_{1}=\|\varphi_{1}-\varphi_{0}\|_{1}=\int_0^{460}|p(x)|dx=0.00008070175.
\eqno (4.27)
$$
then, using (3.9) and (4.6) we can calculate that the corresponding global contraction constant is $$q=0.00008070175\times\frac{1}{1.00787\times10^{-4}}=0.8007.$$With this and the iterates $\varphi_{2},\varphi_{3}$, (4.2) yields the following a posteriori estimate:
$$
\aligned
\|\varphi_2-\varphi\|_1&\leq \frac{1}{1-q}\|\varphi_{3}-\varphi_2\|_1=\frac{1}{1-0.8007}\times2.648099\times10^{-7}=1.3287\times10^{-6}.
\endaligned
$$
Consequently, from (4.3) in Theorem 4.1 we have
$$
\|w_{2}-w\|_{1}\leq K \|\varphi_{2}-\varphi\|_{1}\leq K\frac{1}{1-q}\|\varphi_{3}-\varphi_2\|_1
=3.014\times10^{7}\times1.3287\times10^{-6}=40.04702;
$$Since we can calculate $\|w_2-w_1\|_1\approx 4.570394$, then
$$
\|w-w_1\|_1 \leq \|w_2-w_1\|_1 + \|w-w_2\|_1\leq4.570394+40.04702=44.617414.
$$
Therefore, based on the above estimate obtained from the global contraction constant $q=0.8007$, we can see that the mean absolute errors between the exact solution $w$ and the first two iterates $w_{1}$ and $w_{2}$ over the entire main span do not exceed $\frac{44.617414}{460}=0.09699437m$ and $\frac{40.04702}{460}=0.08705874m$, respectively. Comparing the above results with the amplitudes of $w_{1}$ and $w_{2}$ computed from (4.23) and (4.25):
$$\underset{x\in (0,460)}\max w_{1}(x)=w_{1}(230)=0.46873 m,\ \underset{x\in (0,460)}\min w_{1}(x)=w_{1}(53.58853)=w_{1}(406.41147)=-0.13627 m,\ $$
$$\underset{x\in (0,460)}\max w_{2}(x)=w_{2}(230)=0.44986 m,\ \underset{x\in (0,460)}\min w_{2}(x)=w_{2}(52.81768)=w_{2}(407.18233)=-0.12383 m
$$
we see that the global error bounds are rigorous and still useful, but
somewhat conservative when compared with the displacement amplitudes.  This
is mainly because the global contraction constant is close to \(1\) in the
present case.  To obtain sharper estimates that better reflect the actual
numerical error levels suggested by the iteration, we next use the observed
local contraction factor.  By direct calculation, we have
$$
\frac{\|\varphi_2-\varphi_1\|_1}
{\|\varphi_1-\varphi_0\|_1}
\approx
\frac{4.166122\times10^{-6}}
{8.070175\times10^{-5}}
\approx5.162369\times10^{-2}
=:q_1,
$$
and
$$
\frac{\|\varphi_3-\varphi_2\|_1}
{\|\varphi_2-\varphi_1\|_1}
\approx
\frac{2.648099\times10^{-7}}
{4.166122\times10^{-6}}
\approx6.356269\times10^{-2}
=:q_2.
$$
Therefore, we take
$$
q_{\rm loc}:=\max\{q_1,q_2\}
=6.356269\times10^{-2}
$$
as the observed local contraction factor and use it to obtain the following
a posteriori estimate:$$
\|\varphi-\varphi_2\|_1
\approx
\frac{1}{1-q_{\rm loc}}\|\varphi_3-\varphi_2\|_1
\approx
\frac{2.648099\times10^{-7}}
{1-6.356269\times10^{-2}}
=2.827844\times10^{-7}.
$$
Consequently, we have
$$
\|w-w_2\|_1\leq K\|\varphi-\varphi_2\|_1\approx3.0147\times10^7\times2.827844\times10^{-7}=8.525101.
$$
Combining this with $\|w_2-w_1\|_1\approx 4.570394$, we further estimate
$$
\|w-w_1\|_1 \leq
\|w_2-w_1\|_1+\|w-w_2\|_1\approx4.570394+8.525101=13.095495.
$$
In terms of mean absolute errors, we
obtain
$$
\frac{\|w-w_1\|_1}{460}\leq 2.846847\times10^{-2}m,\qquad
\frac{\|w-w_2\|_1}{460}\leq 1.853283\times10^{-2}m.
\qquad
$$
Thus, compared with the estimates obtained from the global contraction constant,
the observed local factor gives much smaller a posteriori numerical
estimates.  In particular, the mean absolute error estimate for the second
iterate is about \(1.85\times10^{-2}\,{\rm m}\), which is small compared
with the displacement amplitude of \(w_2\).  This indicates that \(w_2\)
provides a reasonable engineering-level numerical approximation for the
present load case.

It should be emphasized that the estimates based on \(q_{\rm loc}\) are
not rigorous upper bounds in the sense of the global contraction theorem,
but they provide a more realistic numerical indication of the actual
numerical local convergence behavior.

\noindent\textbf{Remark 4.1.}
The explicit expressions of \(w_1,w_2\) and the corresponding plots show
that, in both centered live-load cases, the displacement curves are
symmetric with respect to the midspan.  The maximum downward deflection
occurs at the midspan, while the maximum upward deflections and the two
inflection points are located symmetrically on the two sides of the span.

In both cases, the second iterate \(w_2\) has a smaller amplitude than the
first iterate \(w_1\), which suggests that the solution \(w_1\) of the simplified linear nonlocal model tends to overestimate the displacement level of the original classical nonlinear nonlocal model. For the narrow small live load, even the conservative theoretical global error
estimate gives mean errors of about \(1.34\%\) for \(w_1\) and
\(2.08\times10^{-3}\%\) for \(w_2\), measured relative to the maximum
downward displacement amplitude \(\max w_2\).  Thus, \(w_1\) is already a
very accurate approximation.  For the wider and larger live load, the
global theoretical estimate is less sharp because the corresponding global
contraction constant is close to \(1\), namely \(q\approx0.8007\).  However,
the estimates based on the observed local contraction factor show that the
mean errors are about \(6.33\%\) for \(w_1\) and \(4.12\%\) for \(w_2\),
again measured relative to the maximum downward displacement amplitude
\(\max w_2\).  Hence, even in this larger-load case, the first iterate
\(w_1\) still gives an acceptable first approximation to the nonlinear
model.

\noindent\textbf{Remark 4.2.}
The comparison between the global contraction constant and the observed
local contraction factor shows the gap between the worst-case theoretical
estimate and the actual numerical convergence of the iteration.  In the
narrow small live-load case,
\[
q_{\rm loc}=2.59304\times10^{-3}
<q=0.01740687,
\]
while in the wider and larger live-load case,
\[
q_{\rm loc}=6.356269\times10^{-2}
< q=0.8007.
\]
Thus, in both centered live-load cases, the observed local contraction factor is about one order of magnitude smaller than the corresponding global
contraction constant. This shows that the global error estimates, although rigorous and useful, are conservative.  Using \(q_{\rm loc}\) as the effective contraction factor gives smaller a posteriori estimates that better reflect the actual numerical convergence of the iteration, although they are not rigorous upper bounds in the sense of the global contraction theorem.

\noindent{\bf Example 4.2.} We consider a long actual suspension bridge---the famous Golden Gate Bridge, which has a main span of $L=1280$ m with the sag of $d=145$ m. The static behaviour of this bridge has been considered by [14]. In Table 1 of their paper,
the authors reported the main structural parameters of the Golden Gate
Bridge, including the dead load $q=307$ and the two dimensionless parameters: \[
\lambda^2 = \frac{q^2 E_c A_c L^3}{H^3 L_c}=231, \qquad
\alpha^2 = \frac{EI}{H L^2}=0.001523
\]
which are known as the Irvine parameter and the Steinman stiffness factor, respectively (see [8] and [21]).
Combining these data with the expression $H=\frac{qL^{2}}{8d}$ (see (1) in [14]), Melan equation (1.4) can be transformed into
the equivalent equation (3.1) and the corresponding parameters are:
$$
\aligned
&a=\frac{H}{EI}=\frac{1}{\alpha^{2}L^{2}}=0.0004007561,\ \
b=\frac{q}{H}\frac{E_{c}A_{c}}{EI}\frac{1}{L_{c}}=c\frac{H}{q}=c\frac{L^{2}}{8d}=6.234824\times10^{-8},\\ &c=(\frac{q}{H})^2\frac{E_{c}A_{c}}{EI}\frac{1}{L_{c}}=\frac{\lambda^{2}}{\alpha^{2}L^{5}}=4.414304\times10^{-11}.\\
\endaligned \eqno (4.28)
$$
By (H), we can calculate the following:
$$
\frac{2a^2
\left[1+c\left(\frac{L^3}{12a}-\frac{L}{a^2}+\frac{2}{a^{5/2}}\tanh\!\left(\frac{\sqrt aL}{2}\right)\right)\right]^2}{
bL^2\left[
1+
c\left(
\frac{5L^3}{24a}
-\frac{L}{a^2}
+
\frac{1}{a^{5/2}}
\left(
2-\frac{aL^2}{4}
\right)
\tanh\!\left(\frac{\sqrt aL}{2}\right)
\right)
\right]}=2.682132\times10^{-5},
\eqno (4.29)$$
that is, the function $p$ in (3.1) should satisfy $\|p\|_{1}<2.682132\times10^{-5}$
and then the real live load $P$ on the bridge in (1.4) should meet the following condition:
$$
\aligned
\|P\|_{1}&<2.682132\times10^{-5}\times EI=2.682132\times10^{-5}\times \alpha^{2}\frac{L^{4}q}{8d}\\
&=2.682132\times10^{-5}\times 1.081982\times10^{9}=29020.19(\mathrm{kN}).
\endaligned\eqno (4.30)
$$

We now consider three live-load cases for the Golden Gate Bridge.  We first
take \(30\,{\rm kN/m}\) as a representative moderate live-load intensity.
Two typical loading patterns are considered at this load level: an
asymmetric left-half-span loading condition and a symmetric central-half-span
loading condition.  This allows us to compare the influence of the loading
position on the displacement response and on the convergence behaviour of
the iteration.  In addition, motivated by the \(60\,{\rm kN/m}\) load level
considered in the existing literature, we also examine a heavier
left-half-span load of \(60\,{\rm kN/m}\) as a supplementary numerical
comparison.

\noindent\textbf{Case 1.} A uniformly distributed live load of
\(P_0=30\,{\rm kN/m}\) acts on the left half of the main span, namely
\(0\leq x\leq640\).  This load may be interpreted as an equivalent traffic
load corresponding approximately to three loaded traffic lanes, each
carrying \(10\,{\rm kN/m}\), over the left half of the main span. The actual live load function $P(x)$ in (1.4) for this case is:
$$
P(x)=
\left\{
\begin{aligned}
30,\quad &0\leq x\leq640,\\
0,\quad &640<x\leq1280.
\end{aligned}
\right.
\eqno (4.31)
$$
Since the total actual live load is $$
\|P\|_1=\int_0^{1280} |P(x)|\,dx=P_0\frac{1280}{2}=30\times640=19200\ {\rm kN}.
$$ it follows from (4.30) that the Condition (H) of Theorem 3.1 holds in this case.
By $EI=1.081982\times10^{9}$, the function $p(x)$ in equation (3.1) is
$$
p(x)=
\left\{
\begin{aligned}
\dfrac{30}{1.081982\times10^{9}},\quad &0\leq x\leq640,\\
0,\quad &640<x\leq1280.
\end{aligned}
\right.
\eqno (4.32)
$$
For this normalized load and the parameters \(a,b,c\) given in (4.28), the problem corresponding to (3.1) becomes
$$
\aligned
&w''''(x)-\left(a+b\int_0^{1280} w(x)dx\right)w''(x)
+c\int_0^{1280} w(x)dx=p(x),\quad x\in (0,1280),\\
&w(0)=w(1280)=w''(0)=w''(1280)=0.
\endaligned
\eqno (4.33)
$$

Following the same iterative procedure as in Example 4.1, we take $\varphi_{0}\equiv0$ as the initial function and then, by (2.6) and (3.4), successively compute the following:
$$
\aligned
&w_{0}\equiv0;\\
&\varphi_{1}(x)=p(x);\\
\endaligned
$$
\begin{small}
$$
\aligned
&w_{1}(x)
=p\left[\int^{640}_{0}G(x,s)ds
-k\int^{1280}_{0}G(x,s)ds\right]\\
&=p\left\{
\begin{aligned}
&
\frac{\left(\frac34-k\right)1280x+(k-1)x^{2}}{2a}
\\
&+\frac{
\sinh(\sqrt a x)\left(\cosh\frac{\sqrt a1280}{2}-k\right)
+(1-k)\left(\sinh(\sqrt a(1280-x))-\sinh(\sqrt a1280)\right)}
{a^{2}\sinh(\sqrt a1280)},\quad 0\leq x\leq640,\\[1mm]
&\frac{\frac{1280^2}{4}-\left(\frac14+k\right)1280x+kx^{2}}{2a}\\
&-\frac{\sinh(\sqrt a(1280-x))\left(\cosh\frac{\sqrt a1280}{2}-1+k\right)
+k\left(\sinh(\sqrt ax)-\sinh(\sqrt a1280)\right)}{a^{2}\sinh(\sqrt a1280)},
\quad 1280\geq x>640.
\end{aligned}
\right.
\endaligned
\eqno (4.34)
$$
\end{small}
where $p=\frac{30}{1.081982\times10^{9}}, k=
\frac{c\int^{1280}_{0}\int^{640}_{0}G(x,s)dsdx}
{1+c\int^{1280}_{0}\int^{1280}_{0}G(x,s)dsdx}
=0.4749067, a=0.0004007561;$
\begin{small}
$$
\aligned
&\varphi_{2}(x)
=p(x)+b\int_0^L w_{1}(x)dx\,w_{1}''(x)=p(x)+6.234824\times10^{-8}\times298.2958\,w_{1}''(x)\\
&=\left\{
\begin{aligned}
&
\frac{30}{1.081982\times10^{9}}
+0.00001859822\frac{p}{a}
\left[
(k-1)
+\frac{
\sinh(\sqrt ax)\left(\cosh\frac{\sqrt a1280}{2}-k\right)}
{\sinh(\sqrt a1280)}
\right.
\\
&\hspace{2.5em}
\left.
+\frac{
(1-k)\sinh(\sqrt a(1280-x))}
{\sinh(\sqrt a1280)}
\right],
\quad\quad \quad \quad\quad \quad 0\leq x\leq640,
\\[2mm]
&
0.00001859822\frac{p}{a}
\left[
k
-\frac{
\sinh(\sqrt a(1280-x))
\left(\cosh\frac{\sqrt a1280}{2}-1+k\right)
+k\sinh(\sqrt ax)}
{\sinh(\sqrt a1280)}
\right],
\quad640<x\leq1280
\end{aligned}
\right.
\\
&\approx
\left\{
\begin{aligned}
&2.705122\times10^{-8}+3.509702\times10^{-15}\sinh(\sqrt ax)+1.005349\times10^{-20}\sinh(\sqrt a(1280-x)),\quad0\leq x\leq640,\\[1mm]
&6.110835\times10^{-10}-9.092614\times10^{-21}\sinh(\sqrt ax)-3.509701\times10^{-15}\sinh(\sqrt a(1280-x)),\quad640<x\leq1280;
\end{aligned}\right.\endaligned
\eqno (4.35)
$$
\end{small}
\begin{small}
$$
\aligned
&w_{2}(x)
=\int^{L}_{0}G(x,s)\varphi_{2}(s)ds-\frac{c\int^{L}_{0}\int^{L}_{0}G(x,s)\varphi_{2}(s)dsdx}
{1+c\int^{L}_{0}\int^{L}_{0}G(x,s)dsdx}\int^{L}_{0}G(x,s)ds\\
&=\int^{L}_{0}G(x,s)\varphi_{2}(s)ds-1.313752\times10^{-8}\int^{L}_{0}G(x,s)ds\\
&\approx\left\{\begin{aligned}&-0.08242587+0.01166382x-0.00001735932x^{2}\\
&+2.871927\times10^{-7}\sinh(\sqrt a x)+1.226455\times10^{-12}\sinh(\sqrt a(1280-x))\\
&+2.187359\times10^{-10}x\sinh(\sqrt a x)-6.265625\times10^{-16}x\sinh(\sqrt a(1280-x)),\quad 0\leq x\leq640,\\[1mm]
&13.586-0.03056059x+0.0000156285x^2\\
&-3.785654\times10^{-13}\sinh(\sqrt a x)-5.671745\times10^{-7}\sinh(\sqrt a(1280-x))\\
&-5.666776\times10^{-16}x\sinh(\sqrt a x)+2.187358\times10^{-10}x\sinh(\sqrt a(1280-x)),\quad 640<x\leq1280;
\end{aligned}\right.
\endaligned\eqno (4.36)
$$\end{small}
{\small
\begin{equation}
\begin{aligned}
&\varphi_{3}(x)
=p(x)+b\int_{0}^{1280}w_{2}(x)\,dx\,w_{2}''(x)\\
&=p(x)+6.234824\times10^{-8}\times297.6125\,w_{2}''(x)=p(x)+1.855561\times10^{-5}w_{2}''(x)\\
&\approx
\left\{\begin{aligned}
&2.708266\times10^{-8}+2.135644\times10^{-15}\sinh(\sqrt a x)+9.120259\times10^{-21}\sinh(\sqrt a(1280-x))\\
&\quad+1.62658\times10^{-18}x\sinh(\sqrt a x)-4.659292\times10^{-24}x\sinh(\sqrt a(1280-x))\\
&\quad+1.625045\times10^{-16}\cosh(\sqrt a x)+4.654894\times10^{-22}\cosh(\sqrt a(1280-x)),\quad 0\leq x\leq640,\\[2mm]
&5.799929\times10^{-10}-2.815117\times10^{-21}\sinh(\sqrt a x)-4.217666\times10^{-15}\sinh(\sqrt a(1280-x))\\
&\quad-4.213971\times10^{-24}x\sinh(\sqrt a x)+1.62658\times10^{-18}x\sinh(\sqrt a(1280-x))\\
&\quad-4.209993\times10^{-22}\cosh(\sqrt a x)-1.625045\times10^{-16}\cosh(\sqrt a(1280-x)),\quad 640<x\leq1280.
\end{aligned}
\right.
\end{aligned}
\tag{4.37}
\end{equation}
}
Similarly to Example 4.1, we omit the explicit expression for $w_{3}$ owing to the tedious integration involved. The results of the two iterations are shown in Figure 4 below drawn with MATLAB:
\begin{figure}[H]
\centering
\scalebox{1}[0.85]{%
\includegraphics[width=0.9\linewidth]{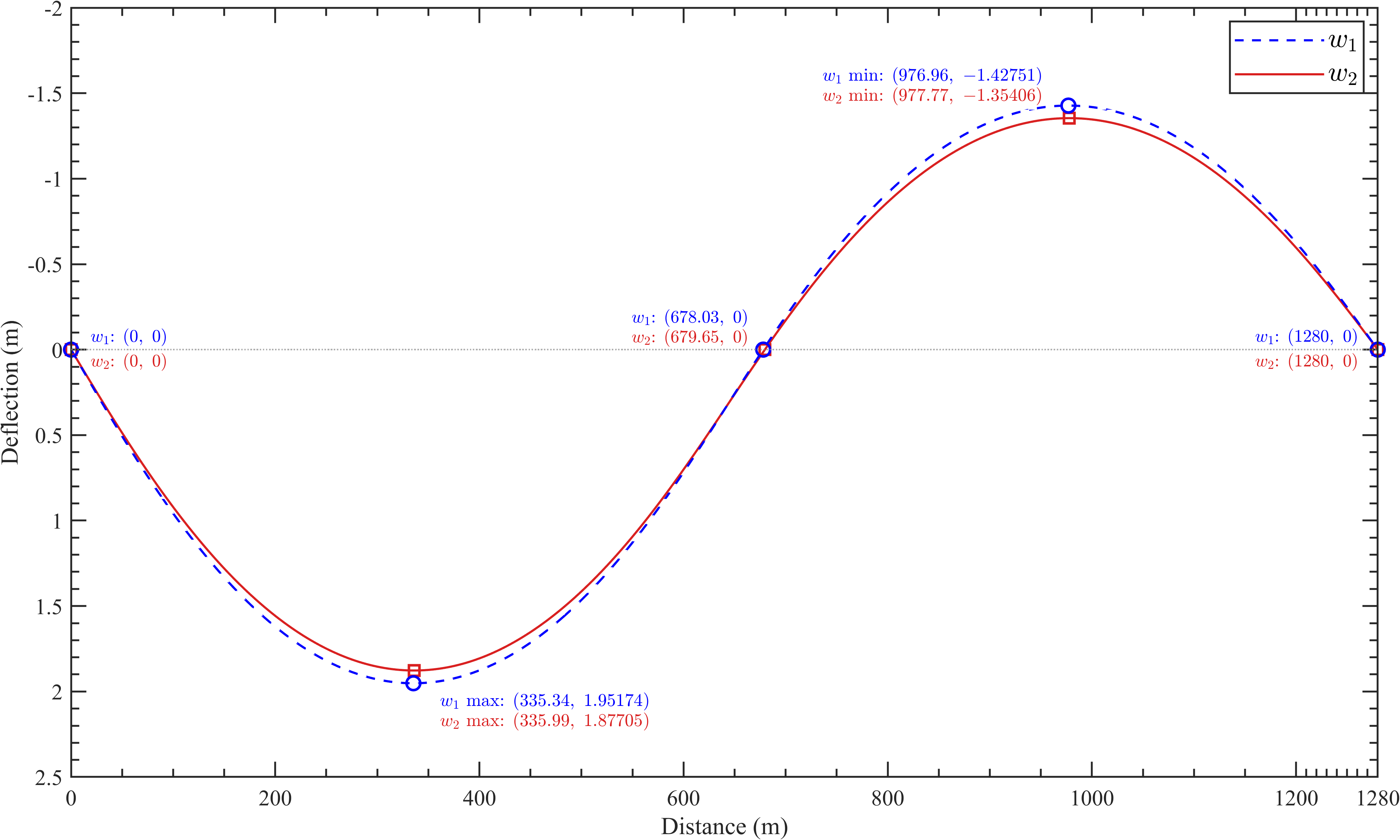}
}
\caption{The first two iterates \(w_1\) and \(w_2\) for the Golden Gate
Bridge under a uniformly distributed live load of \(30\,{\rm kN/m}\) on
the left half of the main span.}
\end{figure}
 As will be seen below, $w_{1}$ already gives an acceptable approximation, while the second iterate $w_{2}$ provides a more accurate approximation. We next quantify the errors of the first two iterates with respect to the exact solution, using both the theoretical global contraction constant and the numerically observed local contraction factor.

According to the load function (4.32), we have
\[
\|p\|_1=\|\varphi_1-\varphi_0\|_1=\int_0^{1280}|p(x)|dx=\frac{30\times640}{1.081982\times10^9}
\approx1.774521\times10^{-5},
\tag{4.38}
\]
then, using (3.9) and (4.29) we can calculate that the corresponding global contraction constant is
\[
q=1.774521\times10^{-5}\times\frac{1}{2.682132\times10^{-5}}=0.6616083.
\tag{4.39}
\]
With this and the iterates $\varphi_{2},\varphi_{3}$, (4.2) yields the following a posteriori estimate:
\[
\begin{aligned}
\|\varphi_2-\varphi\|_1
&\leq\frac{1}{1-q}\|\varphi_3-\varphi_2\|_1=\frac{1}{1-0.6616083}\times3.079405\times10^{-8}\approx9.10012\times10^{-8}.
\end{aligned}
\tag{4.40}
\]
Consequently, combining the above results with the value of $K=9.96421\times10^{8}$ computed from (4.4), Theorem 4.1 yields the following estimate for the iterative solutions of the problem (4.33):
\[
\begin{aligned}
\|w_2-w\|_1
&\leq K\|\varphi_2-\varphi\|_1
\leq9.96421\times10^8\times9.10012\times10^{-8}\approx90.67551.
\end{aligned}
\tag{4.41}
\]
Since we can calculate $\|w_2-w_1\|_1\approx 61.88862$, then
we further obtain
\[
\begin{aligned}
\|w-w_1\|_1
&\leq\|w_2-w_1\|_1+\|w-w_2\|_1\leq61.88862+90.67551=152.56413.
\end{aligned}
\tag{4.42}
\]
Therefore, based on the above estimate obtained from the global contraction constant $q=0.6616083$, we can see that the mean absolute errors between the exact solution $w$ and the first two iterates $w_{1}$ and $w_{2}$ over the entire main span do not exceed $\frac{152.56413}{1280}
=1.191907\times10^{-1}\ {\rm m}$ and $\frac{90.67551}{1280}
=7.084024\times10^{-2}\ {\rm m}$, respectively. Comparing the above results with the amplitudes of $w_{1}$ and $w_{2}$ computed from (4.34) and (4.36):
$$
\underset{x\in(0,1280)}\max w_{1}(x)=w_{1}(335.33834)=1.95174\,{\rm m};
\quad\underset{x\in(0,1280)}\min w_{1}(x)=w_{1}(976.96228)=-1.42751\,{\rm m},
$$
$$
\underset{x\in(0,1280)}\max w_{2}(x)=w_{2}(335.99354)=1.87705\,{\rm m};
\quad\underset{x\in(0,1280)}\min w_{2}(x)=w_{2}(977.76759)=-1.35406\,{\rm m},
$$
we see that the mean absolute error of \(w_1\) is already much smaller than the displacement amplitude, and $w_{2}$ is even more accurate.

In fact, even the acceptable error estimate obtained above from the global contraction constant is still conservative.
 To examine the actual convergence rate and errors, we may define the observed local contraction factor
$$
\aligned
q_{\rm loc}:&=\max\{
\frac{\|\varphi_2-\varphi_1\|_1}{\|\varphi_1-\varphi_0\|_1},\frac{\|\varphi_3-\varphi_2\|_1}{\|\varphi_2-\varphi_1\|_1}\}
=\max\{3.916812\times10^{-2},4.430503\times10^{-2}\}\\
&=4.430503\times10^{-2}
\endaligned\eqno (4.43)
$$
which is about one order of magnitude smaller than the global constant $q=0.6616083.$ Using this local factor, we obtain the following a posteriori estimate:
\[\begin{aligned}
\|\varphi-\varphi_2\|_1&\approx
\frac{1}{1-q_{\rm loc}}\|\varphi_3-\varphi_2\|_1\approx\frac{3.079405\times10^{-8}}{1-4.430503\times10^{-2}}
=3.222163\times10^{-8}.
\end{aligned}
\tag{4.44}
\]
Consequently, we have
\[
\begin{aligned}
\|w-w_2\|_1
&\leq K\|\varphi-\varphi_2\|_1\approx9.96421\times10^8\times3.222163\times10^{-8}=32.10631.
\end{aligned}
\tag{4.45}
\]
Combining this with \(\|w_2-w_1\|_1\approx61.88862\), we further estimate
\[
\begin{aligned}\|w-w_1\|_1&\leq\|w_2-w_1\|_1+\|w-w_2\|_1\leq
61.88862+32.10631=93.99493.
\end{aligned}
\tag{4.46}
\]
In terms of mean absolute errors, we
obtain
\[
\frac{\|w-w_1\|_1}{1280}
\approx7.343354\times10^{-2}\ {\rm m},
\qquad
\frac{\|w-w_2\|_1}{1280}
\approx2.508305\times10^{-2}\ {\rm m}.
\]
Thus, compared with the estimates obtained from the global
contraction constant, the observed local factor gives much smaller
a posteriori numerical estimates.
It should be emphasized that the estimates based on \(q_{\rm loc}\) are
not rigorous upper bounds in the sense of the global contraction theorem,
but they provide a more realistic numerical indication of the actual
local convergence behavior.

\noindent\textbf{Case 2.} A uniformly distributed live load of
\(P_0=30\,{\rm kN/m}\) acts on the central half of the main span, namely
\(320\leq x\leq960\).  This load has the same intensity and the same
loaded length as in Case 1, but is placed symmetrically about the midspan.
Hence the total live load remains unchanged, and
the global contraction condition verified in Case 1 is still satisfied.
This case is used to compare the influence of the loading position on the
displacement response and on the convergence behaviour of the iteration. In this case, the function $p(x)$ in equation (4.33) is
$$
p(x)=
\left\{
\begin{aligned}
0,\quad &0\leq x\leq320,\\
\dfrac{30}{1.081982\times10^9},\quad &320<x\leq960,\\
0,\quad &960<x\leq1280.
\end{aligned}
\right.
\eqno (4.47)
$$
Following the same iterative procedure as in the previous cases, and in
particular as in the central-half-span live-load case of Example 4.1,
where a similar uniformly distributed live-load pattern was considered,
we take $\varphi_{0}\equiv0$ as the initial function and then, by (2.6) and (3.4), successively compute the following:
$$
\aligned
&w_{0}\equiv0;\\
&\varphi_{1}(x)=p(x);\\
\endaligned
$$
\begin{small}
$$
\aligned
&w_{1}(x)
=p\left[
\int^{960}_{320}G(x,s)\,ds
-k\int^{1280}_{0}G(x,s)\,ds
\right]
\\
&=p\left\{
\begin{array}{@{}l@{}}
\displaystyle
\frac{1280x(1-2k)+2kx^{2}}{4a}
+\frac{1}{a^{2}\sinh(1280\sqrt a)}
\bigg[
\sinh(\sqrt a x)
\left(
\cosh(320\sqrt a)-\cosh(960\sqrt a)-k
\right)
\\
\displaystyle
\hspace{2.0em}
-k\sinh(\sqrt a(1280-x))
+k\sinh(1280\sqrt a)
\bigg],
\hspace{3.2em}0\leq x\leq320,
\\[4mm]
\displaystyle
\frac{1280x(1-k)+(k-1)x^{2}-\frac{1280^{2}}{16}}{2a}
+\frac{1}{a^{2}\sinh(1280\sqrt a)}
\bigg[
\left(\sinh(\sqrt a(1280-x))+\sinh(\sqrt a x)\right)
\left(\cosh(320\sqrt a)-k\right)
\\
\displaystyle
\hspace{2.0em}
-(1-k)\sinh(1280\sqrt a)
\bigg],
\hspace{8.3em}320<x\leq960,
\\[4mm]
\displaystyle
\frac{1280^{2}-1280x(1+2k)+2kx^{2}}{4a}
+\frac{1}{a^{2}\sinh(1280\sqrt a)}
\bigg[
\sinh(\sqrt a(1280-x))
\left(
\cosh(320\sqrt a)-\cosh(960\sqrt a)-k
\right)
\\
\displaystyle
\hspace{2.0em}
-k\sinh(\sqrt a x)
+k\sinh(1280\sqrt a)
\bigg],
\hspace{5.8em}960<x\leq1280.
\end{array}
\right.
\endaligned
\eqno (4.48)
$$
\end{small}
where $p=\frac{30}{1.081982\times10^{9}}, k=\frac{c\int_0^{1280}\int_{320}^{960}G(x,s)\,ds\,dx}{1+c\int_0^{1280}\int_0^{1280}G(x,s)\,ds\,dx}
\approx0.6553621, a=0.0004007561;$
\begin{small}
$$
\aligned
&\varphi_{2}(x)
=p(x)+b\int_{0}^{1280}w_{1}(x)\,dx\,w_{1}''(x)
=p(x)+6.234824\times10^{-8}\times411.6425\,w_{1}''(x)
\\
&=
\left\{
\begin{array}{@{}l@{}}
\displaystyle
2.566518\times10^{-5}\frac{p}{a}
\left[
k+
\frac{
\left(
\cosh(320\sqrt a)
-\cosh(960\sqrt a)
-k
\right)\sinh(\sqrt a x)
-k\sinh(\sqrt a(1280-x))}
{\sinh(1280\sqrt a)}
\right],
\\[-1mm]
\hfill\hspace{2em}0\leq x\leq320,
\\[3mm]
\displaystyle
\frac{30}{1.081982\times10^{9}}
+
2.566518\times10^{-5}\frac{p}{a}
\left[
k-1+
\frac{
\left(
\cosh(320\sqrt a)-k
\right)
\left(
\sinh(\sqrt a x)+\sinh(\sqrt a(1280-x))
\right)}
{\sinh(1280\sqrt a)}
\right],
\\[-1mm]
\hfill\hspace{2em}320<x\leq960,
\\[3mm]
\displaystyle
2.566518\times10^{-5}\frac{p}{a}
\left[
k+
\frac{
\left(
\cosh(320\sqrt a)
-\cosh(960\sqrt a)
-k
\right)\sinh(\sqrt a(1280-x))
-k\sinh(\sqrt a x)}
{\sinh(1280\sqrt a)}
\right],
\\[-1mm]
\hfill\hspace{2em}960<x\leq1280
\end{array}
\right.
\\[2mm]
&\approx
\left\{
\begin{array}{@{}l@{\hspace{2.5em}}r@{}}
\displaystyle
1.163715\times10^{-9}
-2.932606\times10^{-12}\sinh(\sqrt a\,x)-1.731549\times10^{-20}\sinh(\sqrt a(1280-x)), \ \ 0\leq x\leq320,
\\[3mm]
\displaystyle
2.711492\times10^{-8}
+7.981674\times10^{-18}\sinh(\sqrt a\,x)+7.981674\times10^{-18}\sinh(\sqrt a(1280-x)),\ \ 320<x\leq960,
\\[3mm]
1.163715\times10^{-9}
-1.731549\times10^{-20}\sinh(\sqrt a\,x)-2.932606\times10^{-12}\sinh(\sqrt a(1280-x)),\ \ 960<x\leq1280;
\end{array}
\right.
\endaligned
\eqno (4.49)
$$
\end{small}
\begin{small}
$$
\aligned
&w_{2}(x)
=\int^{1280}_{0}G(x,s)\varphi_{2}(s)\,ds
-
\frac{
c\int^{1280}_{0}\int^{1280}_{0}G(x,s)\varphi_{2}(s)\,ds\,dx}
{1+c\int^{1280}_{0}\int^{1280}_{0}G(x,s)\,ds\,dx}
\int^{1280}_{0}G(x,s)\,ds
\\
&=\int^{1280}_{0}G(x,s)\varphi_{2}(s)\,ds
-1.811026\times10^{-8}
\int^{1280}_{0}G(x,s)\,ds
\\
&\approx
\left\{
\begin{aligned}
&
0.09827084
-0.006341527x
+2.114322\times10^{-5}x^{2}
\\
&\quad
-1.901157\times10^{-4}\sinh(\sqrt a x)
-1.46222\times10^{-12}\sinh(\sqrt a(1280-x))
\\
&\quad
-1.827694\times10^{-7}x\sinh(\sqrt a x)
+1.076437\times10^{-15}x\sinh(\sqrt a(1280-x)),
 \ 0\leq x\leq320,
\\[2mm]
&
-3.367743
+0.01438026x
-1.123458\times10^{-5}x^{2}
\\
&\quad
+1.994234\times10^{-10}\sinh(\sqrt a x)
+8.36151\times10^{-10}\sinh(\sqrt a(1280-x))
\\
&\quad
+4.974435\times10^{-13}x\sinh(\sqrt a x)
-4.974435\times10^{-13}x\sinh(\sqrt a(1280-x)),\
320<x\leq960,
\\[2mm]
&
26.62216
-0.04778511x
+2.114322\times10^{-5}x^{2}
\\
&\quad
-8.438054\times10^{-14}\sinh(\sqrt a x)
-4.240605\times10^{-4}\sinh(\sqrt a(1280-x))
\\
&\quad
-1.076437\times10^{-15}x\sinh(\sqrt a x)
+1.827694\times10^{-7}x\sinh(\sqrt a(1280-x)),\
960<x\leq1280.
\end{aligned}
\right.
\endaligned
\eqno (4.50)
$$
\end{small}
{\small \begin{equation} \begin{aligned} &\varphi_{3}(x)
=p(x)+b\int_{0}^{1280}w_{2}(x)\,dx\,w_{2}''(x)
=p(x)+6.234824\times10^{-8}\times410.2632\,w_{2}''(x) \\
&=p(x)+2.557919\times10^{-5}w_{2}''(x) \\
&\approx \left\{ \begin{aligned} &
1.081653\times10^{-9}
-1.948879\times10^{-12}\sinh(\sqrt a x)
-1.498924\times10^{-20}\sinh(\sqrt a(1280-x)) \\
&\quad
-1.873572\times10^{-15}x\sinh(\sqrt a x)
+1.103458\times10^{-23}x\sinh(\sqrt a(1280-x)) \\
&\quad
-1.871803\times10^{-13}\cosh(\sqrt a x)
-1.102416\times10^{-21}\cosh(\sqrt a(1280-x)),
\quad 0\leq x\leq320, \\[2mm] &
2.715214\times10^{-8}
+2.044292\times10^{-18}\sinh(\sqrt a x)
+8.571397\times10^{-18}\sinh(\sqrt a(1280-x)) \\
&\quad
+5.099301\times10^{-21}x\sinh(\sqrt a x)
-5.099301\times10^{-21}x\sinh(\sqrt a(1280-x)) \\
&\quad
+5.094488\times10^{-19}\cosh(\sqrt a x)
+5.094488\times10^{-19}\cosh(\sqrt a(1280-x)),
\quad 320<x\leq960, \\[2mm] &
1.081653\times10^{-9}
-8.649862\times10^{-22}\sinh(\sqrt a x)
-4.34705\times10^{-12}\sinh(\sqrt a(1280-x)) \\
&\quad
-1.103458\times10^{-23}x\sinh(\sqrt a x)
+1.873572\times10^{-15}x\sinh(\sqrt a(1280-x)) \\
&\quad
-1.102416\times10^{-21}\cosh(\sqrt a x)
-1.871803\times10^{-13}\cosh(\sqrt a(1280-x)),
\quad 960<x\leq1280. \end{aligned} \right. \end{aligned} \tag{4.51} \end{equation} }
As in Example 4.1, we omit the explicit expression for \(w_3\). The results of the two iterations are shown in Figure 5 below drawn with MATLAB:
 \begin{figure}[H]
\centering
\includegraphics[width=0.9\linewidth]{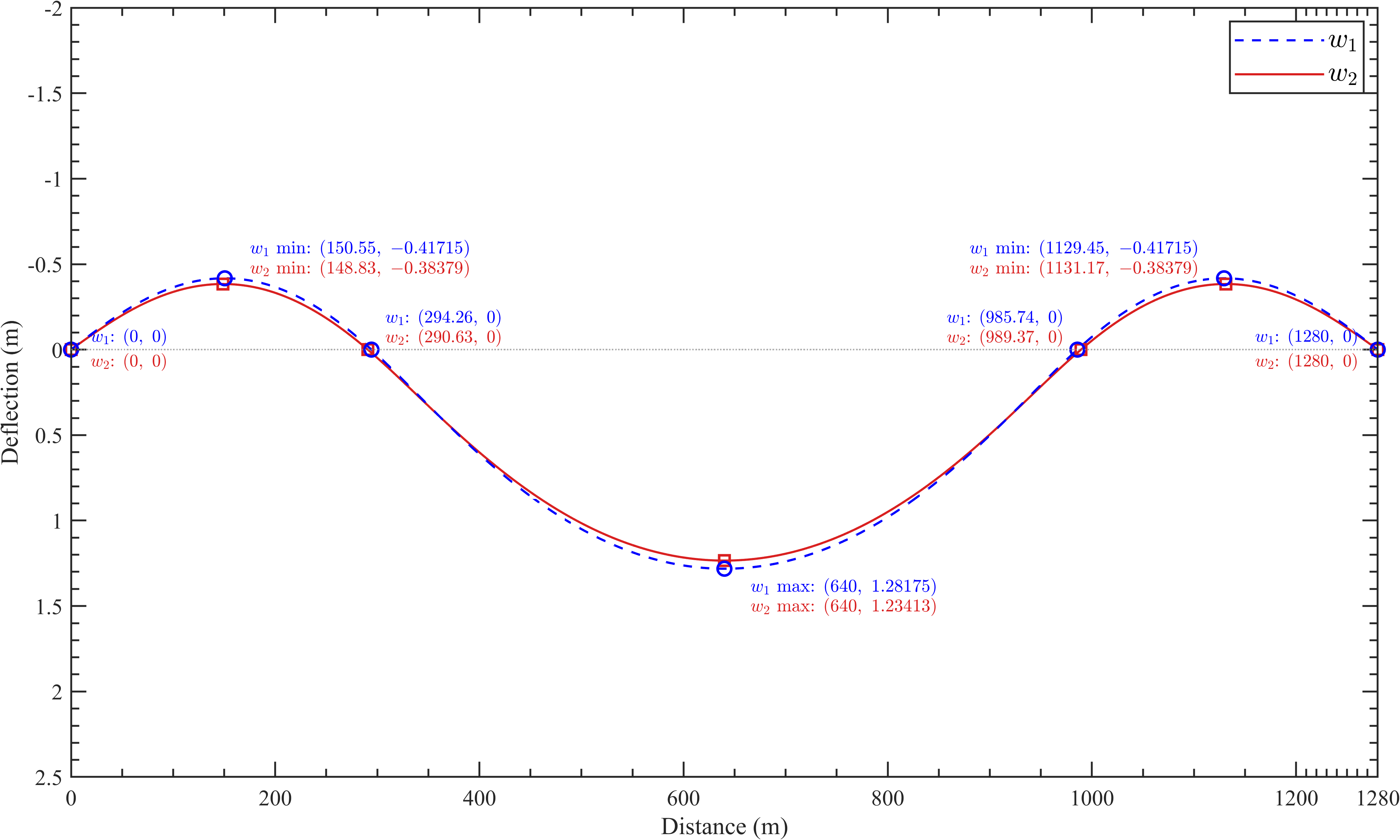}
\caption{The first two iterates \(w_1\) and \(w_2\) for the Golden Gate
Bridge under a uniformly distributed live load of \(30\,{\rm kN/m}\) on
the central half of the main span.}
\end{figure}
We then proceed to estimate the errors of \(w_1\)  and
\(w_2\), using both the global contraction constant and the observed local contraction factor. Since the present central-half-span load has the same intensity and the
same loaded length as in Case 1, we have the same value of the load norm,
namely
\[
\|p\|_1
=
\|\varphi_1-\varphi_0\|_1
=
\int_0^{1280}|p(x)|dx
=
\frac{30\times640}{1.081982\times10^9}
\approx1.774521\times10^{-5}.
\]
Therefore, the global contraction constant is the same as that in Case 1:
\[
q\approx0.6616083.
\]
Combining this value of \(q\) with the iterates \(\varphi_2\) and
\(\varphi_3\) computed in the present case, (4.2) yields the following
a posteriori estimate:
\[
\begin{aligned}
\|\varphi_2-\varphi\|_1\leq\frac{1}{1-q}\|\varphi_3-\varphi_2\|_1=\frac{1}{1-0.6616083}
\times4.869838\times10^{-8}\approx1.439113\times10^{-7}.
\end{aligned}
\tag{4.52}
\]
Since the constant \(K\) defined in (4.4) depends only on the structural parameters
\(a,c,L\), it is the same as in Case 1, namely \(K\approx9.96421\times10^{8}\). Consequently, Theorem 4.1 gives
\[\begin{aligned}\|w_2-w\|_1\leq K\|\varphi_2-\varphi\|_1
\leq9.96421\times10^8
\times1.439113\times10^{-7}\approx143.3962.
\end{aligned}
\tag{4.53}
\]
Moreover,
\[
\|w_2-w_1\|_1\approx33.15275,
\tag{4.54}
\]
and therefore
\[\begin{aligned}\|w-w_1\|_1&\leq\|w_2-w_1\|_1+\|w-w_2\|_1\leq33.15275+143.3962
=176.54895.\end{aligned}\tag{4.55}
\]
Thus, based on the global contraction constant, the mean absolute errors
between the exact solution \(w\) and the first two iterates \(w_1\) and
\(w_2\) over the entire main span do not exceed $\frac{176.54895}{1280}
\approx1.379289\times10^{-1}\ {\rm m}$ and $\frac{143.3962}{1280}
\approx1.120283\times10^{-1}\ {\rm m}$ respectively. Comparing the above results with the amplitudes
of $w_{1}$ and $w_{2}$ computed from (4.48) and (4.50):
$$
\underset{x\in(0,1280)}{\max}\,w_1(x)=w_1(640)\approx1.28175\,{\rm m},\ \ \underset{x\in(0,1280)}{\min}\,w_1(x)
=w_1(150.553)=w_1(1129.447)
\approx-0.41715\,{\rm m},
$$
$$
\underset{x\in(0,1280)}{\max}\,w_2(x)=w_2(640)\approx1.23413\,{\rm m},\ \ \underset{x\in(0,1280)}{\min}\,w_2(x)
=w_2(148.831)=w_2(1131.169)
\approx-0.38379\,{\rm m},
$$
we see that, although the displacement amplitudes in the present
central-half-span case are smaller than those in Case 1, the global
mean-error bounds are larger. As in the previous cases, we next examine
the actual convergence rate and errors by introducing the observed local
contraction factor
$$
\aligned
q_{\rm loc}:&=\max\left\{
\frac{\|\varphi_2-\varphi_1\|_1}{\|\varphi_1-\varphi_0\|_1},
\frac{\|\varphi_3-\varphi_2\|_1}{\|\varphi_2-\varphi_1\|_1}
\right\}
=\max\{4.806673\times10^{-2},
5.709379\times10^{-2}\}
\\
&=5.709379\times10^{-2}.
\endaligned\eqno (4.56)
$$
This is much smaller than the global contraction constant
\(q=0.6616083\).  Using this local factor, we obtain
$$
\begin{aligned}
\|\varphi-\varphi_2\|_1
&\approx
\frac{1}{1-q_{\rm loc}}\|\varphi_3-\varphi_2\|_1
\approx
\frac{4.869838\times10^{-8}}
{1-5.709379\times10^{-2}}
=5.164711\times10^{-8}.
\end{aligned}
\eqno (4.57)
$$
Consequently, we have
$$
\begin{aligned}
\|w-w_2\|_1&\leq K\|\varphi-\varphi_2\|_1\approx9.96421\times10^8\times5.164711\times10^{-8}=51.46226.
\end{aligned}
\eqno (4.58)
$$
Combining this with \(\|w_2-w_1\|_1\approx33.15275\), we further estimate
$$
\begin{aligned}
\|w-w_1\|_1&\approx\|w_2-w_1\|_1+\|w-w_2\|_1\approx33.15275+51.46226=84.61501.
\end{aligned}
\eqno (4.59)
$$
In terms of mean absolute errors, we obtain
$$
\frac{\|w-w_1\|_1}{1280}
\approx6.610548\times10^{-2}\ {\rm m},
\qquad
\frac{\|w-w_2\|_1}{1280}
\approx4.020489\times10^{-2}\ {\rm m}.
$$
Thus, the observed local factor gives much smaller a posteriori numerical
estimates than those obtained from the global contraction constant,
indicating that the global estimate is rather conservative in the present
case.  It should be emphasized that the estimates based on \(q_{\rm loc}\)
are not rigorous upper bounds in the sense of the global contraction
theorem, but they provide a more realistic numerical indication of the
error levels and the local convergence behavior.

 \noindent\textbf{Case 3.}
A uniformly distributed live load of $P_{0}=60\ \mathrm{kN/m}$ acts on the left half of the main span, namely $0\leq x\leq640$. Since the Golden Gate Bridge carries six traffic lanes, this load may be interpreted as an idealized situation in which all six lanes over half the main span are densely occupied by heavy trucks, with an equivalent load of $10\ \mathrm{kN/m}$ per lane. (In fact, such continuous dense loading cannot realistically occur under ordinary traffic conditions. Even without allowing for actual vehicle spacing, application of the multiple-presence factor $m=0.65$ for more than three loaded lanes [48, Table 3.6.1.1.2-1] gives the conservative upper estimate $0.65\times60\times640=24960\ \mathrm{kN}<29020.19\ \mathrm{kN}$ for the corresponding service-equivalent total live load, while actual vehicle spacing would further reduce it.) The unreduced $60\ \mathrm{kN/m}$ load is therefore deliberately retained as an extreme live-load benchmark rather than as an ordinary service traffic load. It has the same left-half-span loading configuration as that considered in [14] for the Golden Gate Bridge and enables comparison with the two preceding $30\ \mathrm{kN/m}$ moderate load cases.

The dimensional unreduced live-load function $P(x)$ adopted in (1.4) for this benchmark case is
$$
P(x)=\left\{
\begin{aligned}
60,\quad &0\leq x\leq640,\\
0,\quad &640<x\leq1280,
\end{aligned}
\right.
\eqno (4.60)
$$
and the corresponding normalized live-load function $p(x)$ in (3.1) is
$$
p(x)=\left\{
\begin{aligned}
\dfrac{60}{1.081982\times10^{9}},\quad &0\leq x\leq640,\\
0,\quad &640<x\leq1280.
\end{aligned}
\right.
\eqno (4.61)
$$
Therefore, the total unreduced live load adopted in this benchmark calculation and the corresponding normalized live-load norm are given by
$$
\|P\|_1
=\int_0^{1280}|P(x)|dx
=60\times640
=38400\ {\rm kN},
\qquad
\|p\|_1
=
\frac{60\times640}{1.081982\times10^9}
\approx3.549041\times10^{-5}.
$$
These two quantities are twice the corresponding values in Case 1, and they exceed the admissible bounds in (4.30) and (4.29), respectively, since $38400>29020.19$ and $3.549041\times10^{-5}>2.682132\times10^{-5}$. Hence, the key hypothesis (H) of Theorem 3.1, which is also required in Theorem 4.1, is not satisfied for this unreduced benchmark load; consequently, the existence and uniqueness conclusions and the rigorous global error estimates established under (H) are not guaranteed in this case. Indeed, the quantity $q$ defined in (3.9) is
$$
q=
\frac{3.549041\times10^{-5}}
{2.682132\times10^{-5}}
\approx1.323217>1.
$$
Nevertheless, we still apply the same iterative scheme numerically to examine whether it remains effective under this heavier loading condition, which lies beyond the range covered by the global contraction argument. For this purpose, we compute the first two iterates and compare their displacement levels with the value reported in [14] to assess their consistency with the published engineering result.

Following the same iterative procedure as in the previous cases,
we take $\varphi_{0}\equiv0$ as the initial function and then
$$
\aligned
&w_{0}\equiv0;\\
&\varphi_{1}(x)=p(x).\\
\endaligned
$$
Since the present live-load intensity is twice that in Case 1, the first
linear iterate \(w_1\) is simply twice the corresponding first iterate in
Case 1, that is $$w_{1}(x)
=p\left[\int^{640}_{0}G(x,s)ds
-k\int^{1280}_{0}G(x,s)ds\right],$$
where $p=\dfrac{60}{1.081982\times10^{9}}.$ The subsequent nonlinear correction is then computed from the same
recursive formula as before. We obtain
\begin{small}
$$
\aligned
&\varphi_{2}(x)
=p(x)+b\int_0^L w_{1}(x)dx\,w_{1}''(x)
=p(x)+6.234824\times10^{-8}\times596.5916\,w_{1}''(x)
\\
&=\left\{
\begin{aligned}
&
\frac{60}{1.081982\times10^{9}}
+0.00003719644\frac{p}{a}
\left[
(k-1)
+\frac{
\sinh(\sqrt ax)\left(\cosh\frac{\sqrt a1280}{2}-k\right)}
{\sinh(\sqrt a1280)}
\right.
\\
&\hspace{2.5em}
\left.
+\frac{
(1-k)\sinh(\sqrt a(1280-x))}
{\sinh(\sqrt a1280)}
\right],
\quad\quad \quad \quad\quad \quad 0\leq x\leq640,
\\[2mm]
&
0.00003719644\frac{p}{a}
\left[
k
-\frac{
\sinh(\sqrt a(1280-x))
\left(\cosh\frac{\sqrt a1280}{2}-1+k\right)
+k\sinh(\sqrt ax)}
{\sinh(\sqrt a1280)}
\right],
\quad 640<x\leq1280
\end{aligned}
\right.
\\
&\approx
\left\{
\begin{aligned}
&5.275113\times10^{-8}
+1.403881\times10^{-14}\sinh(\sqrt ax)+4.021397\times10^{-20}\sinh(\sqrt a(1280-x)),
\quad  0\leq x\leq640,
\\[1mm]
&2.444334\times10^{-9}
-3.637046\times10^{-20}\sinh(\sqrt ax)-1.40388\times10^{-14}\sinh(\sqrt a(1280-x)),
\quad 640<x\leq1280.
\end{aligned}
\right.
\endaligned
\eqno (4.62)
$$
\end{small}
\begin{small}
$$
\aligned
&w_{2}(x)
=\int^{L}_{0}G(x,s)\varphi_{2}(s)ds
-\frac{c\int^{L}_{0}\int^{L}_{0}G(x,s)\varphi_{2}(s)dsdx}
{1+c\int^{L}_{0}\int^{L}_{0}G(x,s)dsdx}
\int^{L}_{0}G(x,s)ds
\\
&=\int^{L}_{0}G(x,s)\varphi_{2}(s)ds
-2.621471\times10^{-8}\int^{L}_{0}G(x,s)ds
\\
&\approx
\left\{
\begin{aligned}
&
-0.1483995
+0.02229341x
-0.00003310794x^{2}
\\
&\quad
+2.069925\times10^{-7}\sinh(\sqrt a x)
+2.20811\times10^{-12}\sinh(\sqrt a(1280-x))
\\
&\quad
+8.749435\times10^{-10}x\sinh(\sqrt a x)
-2.50625\times10^{-15}x\sinh(\sqrt a(1280-x)),
\quad 0\leq x\leq640,
\\[1mm]
&
25.84127
-0.05804559x
+0.00002965691x^2
\\
&\quad
+9.256115\times10^{-13}\sinh(\sqrt a x)
-1.32692\times10^{-6}\sinh(\sqrt a(1280-x))
\\
&\quad
-2.26671\times10^{-15}x\sinh(\sqrt a x)
+8.749435\times10^{-10}x\sinh(\sqrt a(1280-x)),
\quad 640<x\leq1280;
\end{aligned}
\right.
\endaligned
\eqno (4.63)
$$
\end{small}
{\small
\begin{equation}
\begin{aligned}
&\varphi_{3}(x)
=p(x)+b\int_{0}^{1280}w_{2}(x)\,dx\,w_{2}''(x)
\\
&=p(x)
+6.234824\times10^{-8}
\times593.8583\,w_{2}''(x)
=p(x)+3.702602\times10^{-5}w_{2}''(x)
\\
&\approx
\left\{
\begin{aligned}
&
5.300206\times10^{-8}
+3.071438\times10^{-15}\sinh(\sqrt a x)
+3.276482\times10^{-20}\sinh(\sqrt a(1280-x))
\\
&\quad
+1.298276\times10^{-17}x\sinh(\sqrt a x)
-3.718875\times10^{-23}x\sinh(\sqrt a(1280-x))
\\
&\quad
+1.297051\times10^{-15}\cosh(\sqrt a x)
+3.715365\times10^{-21}\cosh(\sqrt a(1280-x)),
\quad 0\leq x\leq640,
\\[2mm]
&
2.196154\times10^{-9}
+1.37346\times10^{-20}\sinh(\sqrt a x)
-1.968937\times10^{-14}\sinh(\sqrt a(1280-x))
\\
&\quad
-3.363436\times10^{-23}x\sinh(\sqrt a x)
+1.298276\times10^{-17}x\sinh(\sqrt a(1280-x))
\\
&\quad
-3.360262\times10^{-21}\cosh(\sqrt a x)
-1.297051\times10^{-15}\cosh(\sqrt a(1280-x)),
\quad 640<x\leq1280.
\end{aligned}
\right.
\end{aligned}
\tag{4.64}
\end{equation}
}
The results of the two iterations are shown in Figure 6 below drawn with MATLAB:
 \begin{figure}[H]
\centering
\includegraphics[width=0.9\linewidth]{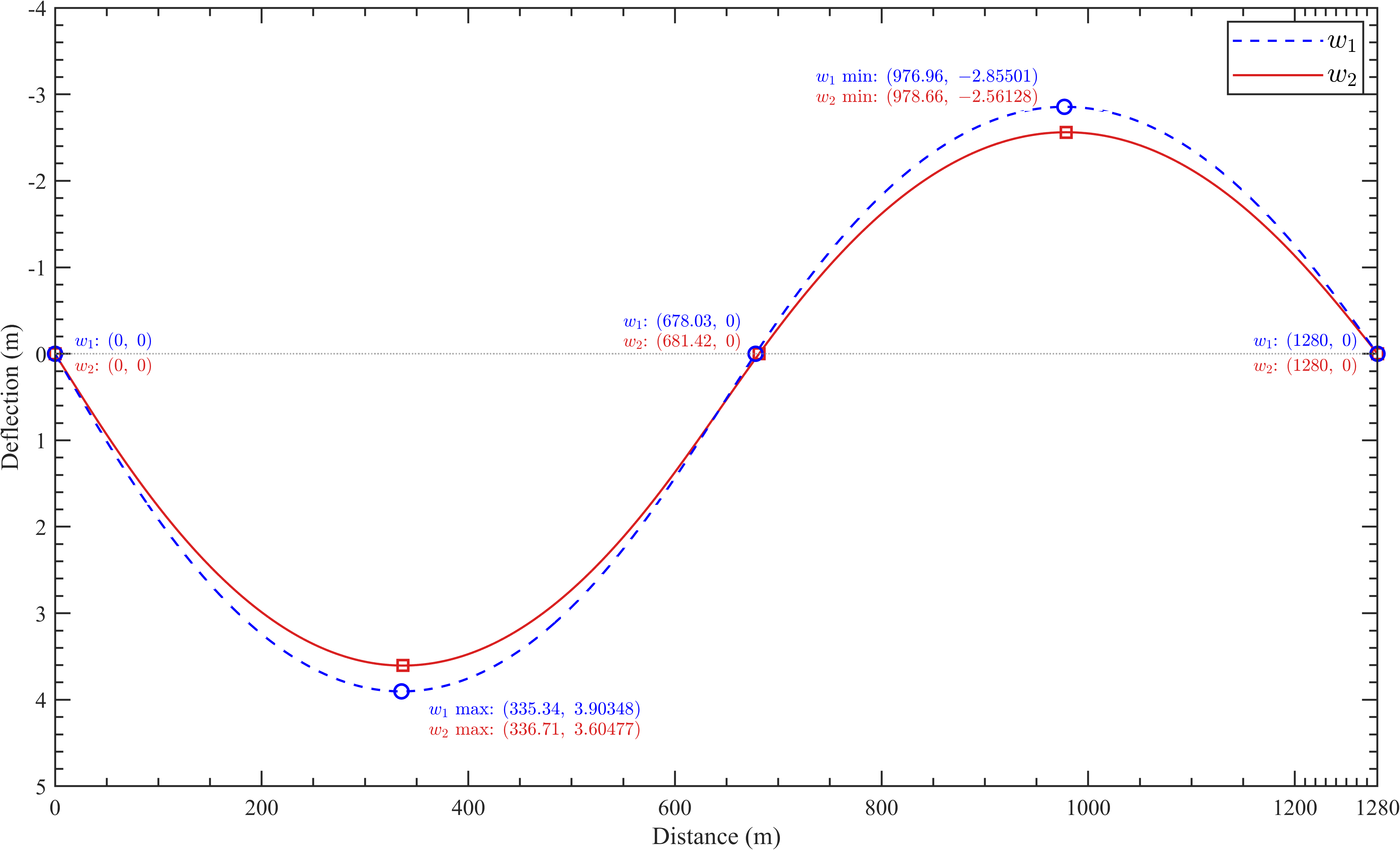}
\caption{The first two iterates \(w_1\) and \(w_2\) for the Golden Gate Bridge under
a uniformly distributed live load of \(60\,{\rm kN/m}\) on the left half of
the main span.}
\end{figure}
From the explicit expressions of \(w_1\) and \(w_2\), we obtain
\[
\max_{x\in(0,1280)} w_1(x)
=w_1(335.33834)
\approx3.90348\ {\rm m},
\qquad
\max_{x\in(0,1280)} w_2(x)
=w_2(336.7117)
\approx3.60477\ {\rm m}.
\]
These maximum downward displacements occur near the quarter point of the
main span.  For comparison, [14] reports
\(W_{\max}\approx3.81\,{\rm m}\) for the Golden Gate Bridge under the same \(60\,{\rm kN/m}\) half-span live load. Thus, the maximum displacement levels obtained from the
first two iterates are of the same order as, and reasonably close to, the
engineering reference values.  This comparison suggests that, although the
sufficient condition of the global contraction theorem is not satisfied,
the iteration may still provide meaningful numerical approximations for
this heavier load case.

We next examine the successive differences of the computed iterates
$\varphi_i(i=0,1,2,3)$ to examine the numerical contraction behaviour suggested by the actual iteration.
The observed ratios are
$$
\frac{\|\varphi_2-\varphi_1\|_1}
{\|\varphi_1-\varphi_0\|_1}
\approx
\frac{2.780189\times10^{-6}}
{3.549041\times10^{-5}}
\approx7.833633\times10^{-2}
=:q_1,
$$
and
$$
\frac{\|\varphi_3-\varphi_2\|_1}
{\|\varphi_2-\varphi_1\|_1}
\approx
\frac{2.45816\times10^{-7}}
{2.780189\times10^{-6}}
\approx8.841701\times10^{-2}
=:q_2.
$$
These two observed ratios are much smaller than $1$, indicating a clear
contraction behaviour along the present numerical iteration.  Therefore,
although the sufficient condition for a global contraction is not satisfied
and no global contraction can be concluded theoretically, it is still
reasonable to introduce the observed local contraction factor
$$
q_{\rm loc}:=\max\{q_1,q_2\}
\approx8.841701\times10^{-2}.
\eqno (4.65)
$$
These small observed ratios indicate a rapid decay of the successive corrections along the computed iteration. However, since hypothesis (H) is not satisfied, $q_{\rm loc}$ is only an observed numerical ratio and cannot be used in (4.2) to obtain a rigorous error estimate with respect to an exact solution. Together with the agreement of the computed displacement levels with the reference values in [14], the present results provide numerical evidence for the effectiveness of the first two iterates under this heavier benchmark load, but do not constitute a proof of convergence.

\noindent\textbf{Remark 4.3.}
For the two left-half-span loading cases
considered for the Golden Gate Bridge, namely Case 1 and Case 3, the
explicit expressions of the iterates and the corresponding plots show that
the maximum downward deflection occurs slightly to the right of the left
quarter point of the main span, while the maximum upward deflection occurs
near the three-quarter point on the unloaded side.  An inflection point is
also located between the downward and upward deflection regions.  These
features agree with the deformation pattern reported in [14]
for the Golden Gate Bridge under nonsymmetric half-span loading.

Similar to the preceding short-span suspension bridge example, in all the
three load cases considered above, the second iterate \(w_2\) has a
smaller amplitude than the first iterate \(w_1\).  This suggests that the solution \(w_1\) of the simplified linear nonlocal model tends to overestimate the displacement level of the original classical nonlinear nonlocal model, whereas \(w_2\) gives a corrected approximation
with a smaller deflection amplitude. Using the a posteriori mean-error estimates based on the observed local contraction factors, and measuring them relative to the maximum downward displacement amplitude \(\max w_2\), we obtain about \(3.91\%\) for \(w_1\) and \(1.34\%\) for \(w_2\) in Case 1, about \(5.36\%\) for \(w_1\) and \(3.26\%\) for \(w_2\) in Case 2. For Case 3, the numerical result agrees well with the engineering reference. Thus, in all three Golden Gate Bridge load cases, the first iterate \(w_1\) already gives an acceptable first approximation to the nonlinear model, while the second iterate gives a further improvement.

\noindent\textbf{Remark 4.4.}
The comparison between Case 1 and Case 2 shows the influence of the loading
position.  In these two cases, both the live-load intensity and the loaded
length are the same, namely \(30\,{\rm kN/m}\) acting over one half of the
main span.  However, the left-half-span load in Case 1 produces a larger
maximum downward deflection than the centered half-span load in Case 2.  For
the more accurate second iterate \(w_2\),
\[
\max w_2\approx1.87705\,{\rm m}\quad\hbox{in Case 1},\qquad
\max w_2\approx1.23413\,{\rm m}\quad\hbox{in Case 2}.
\]
This is consistent with the observation in [14] that, for
nonsymmetric loading lengths around \(50\%\)--\(60\%\) of the span, the
left-side live load produces a larger downward displacement than the
centered load, and the maximum deflection occurs near the quarter point. The heavier \(60\,{\rm kN/m}\) left-half-span load in Case 3 gives
\[
\max w_1\approx3.90348\,{\rm m},\qquad
\max w_2\approx3.60477\,{\rm m},
\]
which is close to the reference maximum displacement level reported in [14].

\noindent\textbf{Remark 4.5.}
For a long-span suspension bridge such as the Golden Gate Bridge, the
observed local contraction factors are still much smaller than the
corresponding global factors, as in the preceding short-span example. In
Case 1 and Case 2, although the global contraction constants are both about
\(0.6616\), the observed local factors are only
\[
q_{\rm loc}\approx4.4305\times10^{-2}
\quad\hbox{and}\quad
q_{\rm loc}\approx5.7094\times10^{-2},
\]
respectively, which are smaller by about one order of magnitude.  Thus, the
actual numerical convergence is much faster than that indicated by the
worst-case global estimate.

This phenomenon is even more significant in Case 3.  For the very heavy
\(60\,{\rm kN/m}\) left-half-span load, the sufficient condition (H) for a
global contraction is no longer satisfied, and the corresponding global
factor is greater than $1$.  Therefore, the theorem does not guarantee a
contraction in this case.  Nevertheless, the observed local factor is still
small, namely $q_{\rm loc}\approx8.8417\times10^{-2}.$
Thus, even beyond the range covered by the global contraction theorem, the
computed iteration still exhibits fast local convergence.  Together with
the agreement of the computed displacement level with the reference values
in [14], this supports the numerical effectiveness of the first
two iterates for the present heavier load case.

\section{Conclusions and open problems}

In this paper, by analyzing the analytical solution and global properties
of the solution to the simplified ``less stiff'' model (1.5),(1.2), we derived a priori
estimates for solutions with nonnegative total displacement and established
an existence and uniqueness result for the classical Melan equation
(1.4),(1.2) under an explicit and concise condition. This result extends
the existence and uniqueness results for small solutions obtained by
Gazzola et al. in [38] for the classical Melan problem (1.4),(1.2) under
small span length \(L\) and live load \(P\), as well as the
result of Wang [40] for suspension bridges with short span length \(L\) and
high deck stiffness \(EI\). The examples of actual suspension bridges
further show that, for practical bridges of different span lengths under
realistic small and moderate downward traffic live loads, the key condition
is satisfied, and hence the corresponding deflection curve is uniquely
determined in the engineering relevant class. To the best of our knowledge, after more than one century of the development
of deflection theory, this is the first rigorous proof of the well-posedness
of the classical Melan equation (1.4),(1.2) in this engineering sense. It
also gives a substantial partial answer to the well-posedness questions
raised by Gazzola et al. in [38] and [39].

We also constructed an efficient iterative approximation method by taking
the solution of the simplified ``less stiff'' model as the first iterate.
For realistic suspension bridges of different span lengths under ordinary
traffic live loads, this provides a rigorously convergent algorithm in the
engineering setting. To the best of our knowledge, this is also the first approximation algorithm with a convergence proof for practical suspension bridges of common span ranges under ordinary live loads. It extends the classical engineering approximation methods which generally
lack rigorous analysis and convergence proofs, and also goes beyond the
monotone iterative method in Wang~[40], whose applicable parameter range is
essentially restricted to short span length \(L\) and high deck stiffness
\(EI\). Moreover, it gives a substantial partial answer to the
open problems on the approximability of solutions raised by Gazzola et al.
At the same time, the method establishes a connection between the classical
Melan model (1.4),(1.2) and the simplified ``less stiff'' model (1.5),(1.2) at both the
theoretical and computational levels. It gives a quantitative analysis of
the difference between the two models and clarifies the influence of the
nonlinear nonlocal term in the original classical equation (1.4). The numerical examples show that the first iterate, namely the solution of the simplified model, already gives a practically acceptable approximation to the deflection curve of the original nonlinear model, although it tends to overestimate the displacement level; the second iterate substantially improves the accuracy. %In the tested cases, two iterations are sufficient for engineering accuracy.
This provides a rigorous
theoretical explanation and quantitative error estimates for the numerical
observations of [28] that the simplified
``less stiff'' model appeared to be able to serve as an effective
approximation to the classical Melan model.

Demonstrative examples on two actual suspension bridges with different
span lengths and under different live-load cases were also carried out to
show the applicability and efficiency of the proposed iterative method. The
computed influence of the live-load magnitude, loading interval and loading
position on the deflection curve is consistent with engineering experience
and observations. More importantly, the numerical results indicate that the
proposed iterative method may remain effective even for very heavy live
loads that violate the key condition (H). This observation naturally leads
to further open problems, for example: for very heavy live loads beyond
condition (H), does the classical Melan equation (1.4),(1.2) still admit existence and uniqueness of solutions, and in particular, is the engineering deflection curve still uniquely determined under such heavy live loads? Can the iterative method developed in this paper with the solution of the simplified
``less stiff'' model (1.5),(1.2) as the first iterate still be proved to converge and remain effective under a rigorous analysis? These questions deserve further
investigation.

%In summary, the present work establishes a rigorous connection between the classical Melan equation and the simplified ``less stiff'' approximation in suspension-bridge deflection theory. It proves the well-posedness of the classical model in the engineering sense, provides a rapidly convergent iterative approximation method with error estimates, and quantitatively explains why the simplified model can serve as an effective first approximation. These results strengthen the analytical foundation of deflection theory for practical suspension bridges and open the way to further studies on heavier live loads, sharper load conditions, and more general deflection behaviors.

\noindent{\bf Declaration of Interest Statement}

\noindent The author declares that there are no competing interests, financial or otherwise, that could
influence or bias the work reported in this manuscript. The author confirms sole responsibility
for the conception, analysis, and interpretation of the results presented in this study.

\noindent{\bf Acknowledgements}

\noindent This work was supported by the NSFC(No.12261057,\ No.12461035).

\vskip 12mm

\centerline {\bf REFERENCES}\vskip5mm\baselineskip 0.45cm
\begin{description}
\baselineskip 15pt

\item{[1]}~ J. Melan, Theory of arches and suspension bridges. Myron Clark Publ. Comp., London, 1913 (German original third edition: Handbuch der Ingenieurwis-senschaften 2, 1906).

\item{[2]}~ C. L. Navier, M\'{e}moire sur les ponts suspendus. Imprimerie Royale, Paris, 1823.

\item{[3]}~ W. J. M. Rankine, A manual of applied mechanics. Charles Griffin and Company, London, 1858.

\item{[4]}~ S. G. Buonopane, D. P. Billington, Theory and history of suspension bridge design from 1823 to 1940, Jour. Strut. Engr. 119 (3) (1993), 954-977.

\item{[5]}~ A. Pugsley, The theory of suspension bridge. Edward Arnold Pub. Ltd, London, 1968. %pp. 75-92.

\item{[6]}~ D. M. Brotton, A general computer program for the solution of suspension bridge problems, Structural engineering. 44 (1966), 161-167.

\item{[7]}~ S. A. Saafan, Theoretical analysis of suspension bridges, Journal of the structural division. 92 (1966), 1-12.

%\item{[9]}~ D. J. Peery, An Influence-Line analysis for suspension bridges, Transactions of the American Society of Civil Engineers. , 121 (1956), 463-480.

%\item{[10]}~ A. Jennings, J. E. Mairs, Static analysis of suspension bridges, Journal of the Structural Division. 92 (1966), 337-342.
\item{[8]}~ M. Irvine, Cable structures. MIT University Press, Cambridge, Mass., 1981. (Reprinted, 1992.)
%Research on the stability of long span self-anchored suspension bridges 作者：Hu, Jian-Hua1, 2;Wang, Lian-Hua2;Shen, Rui-Li3;Xiang, Jian-Jun1;Tang, Mao-Lin3 刊名：Journal of Hunan University (Natural Science) ISSN：1000-2472 出版日期：2008 卷号：Vol.35 期号：No.7 页码：12-15

\item{[9]}~ C. C. Ulstrup. Rating and preliminary analysis of suspension
bridges, Journal of Bridge Engineering, ASCE. 119 (1993), 2653-2679.

\item{[10]}~ G. H. Li, Stability and vibration of bridge. 2nd Ed., China Railway Publishing House, Beijing, 1996. (in Chinese)

\item{[11]}~ P. Walter, Cable-suspended bridges. In: Structural steel designer's handbook. New York: McGraw-Hill Inc, 1999, 15.60-15.67.

\item{[12]}~ P. Clemente, G. Nicolosi, A. Raithel, Preliminary design of very long-span suspension bridges, Engineering Structures. 22 (2000), 1699-1706.

\item{[13]}~ G. P. Wollmann, Preliminary analysis of suspension bridges, J. Bridge Eng. 6 (2001), 227-233.

\item{[14]}~ D. Cobo del Arco, A. C. Aparicio, Preliminary static analysis of suspension bridges, Eng. Struct. 23 (2001), 1096-1103.

\item{[15]}~ J. Q. Lei, M. Z. Zheng, G. Y. Xu, Design of suspension bridge. China Communications Press, Beijing, 2002. (in Chinese)

\item{[16]}~ Y. G. Tan, Z. Zhang, Y. Xu, X. Zhao, Preliminary static analysis of self-anchored suspension bridges, Structural engineering and mechanics, 34 (2010), 281-284.

\item{[17]}~ D. H. Choi, S. G. Gwon, H. Yoo, Nonlinear static analysis of continuous multi-span suspension bridges, International Journal of Steel Structures. 13 (2013), 103-115.

\item{[18]}~ S. U. Shin, M. R. Jung, J. Park, M. Y. Kim, A deflection theory and its validation of earth-anchored suspension bridges under live loads, KSCE Journal of Civil Engineering, 19(1) (2015), 200-212.

\item{[19]}~ M. R. Jung, S. U. Shin, M. M. Attard, M. Y. Kim, Deflection theory for self-anchored suspension bridges under Live Load, Journal of Bridge Engineering, 20 (2015): 04014093.

\item{[20]}~ T. V. Karman, M. A. Biot, Mathematical methods in engineering: an introduction to the mathematical treatment of engineering problems. McGraw-Hill, New York, 1940.

\item{[21]}~ D. B. Steinman, A practical treatise on suspension bridges. 2nd Ed., John Wiley, New York, 1929. (Reprinted, 1953.)

\item{[22]}~ M. T. Godard, Recherches sur le Calcul de la resistance des tabliers des Ponts suspendus. Ann. Ponst et Chaussees, 8 (1894), 105-189.

\item{[23]}~ D. J. Peery, An influence-line analysis for suspension bridge. Transactions of ASCE, 121 (1956), 463-480.

\item{[24]}~ H. Bleich, Die Berechnung Verankerter Hangebr$\ddot{u}$cken,  Berlin, 1935.

\item{[25]}~ S. P. Timoshenko, D. H. Young, Theory of structures. McGraw Hill, New York, 1968.

\item{[26]}~ A. Jennings, Gravity stiffness of classical suspension bridges, Journal of Structural Engineering. 109 (1983), 16-36.

%\item{[27]}~ H. H. West, A. R. Robinson, Continuous method of suspension bridge analysis, Journal of the Structural Division. 49 (1968).
\item{[27]}~ H. Ohshima, K. Sato, N. Watanabe, Structural analysis of suspension bridges, Jour. Engr. Mech. 110 (3) (1984), 392-404.

%\item{[9]}~ H. Neukrich, Berechnung der H$\ddot{a}$ngebr$\ddot{u}$cke bei Ber$\ddot{u}$cksichtigung der Verformung des Kabels, Ingenieur-Archiv, 7, 1936.

%\item{[15]}~ S. U. Shin, M. R. Jung, J. Park, et al, A deflection theory and its validation of earth-anchored suspension bridges under live loads. Ksce Journal of Civil Engineering, 2015, 19(1):200-212.

%\item{[18]}~ D. B. Steinman, A practical treatise on suspension bridges: Their design, construction and erection. Wiley, New York, 1953.

%\item{[19]}~ S. P. Timoshenko, The Stiffness of suspension bridges. Trans. ASCE, Vol 94, 1930.

%\item{[20]}~ M. T. Godard, Recherches sur le Calcul de la resistance des tabliers des Ponts suspendus. Ann. Panst et Chaussees, 1894.

%\item{[21]}~ H. Bleich, Die Berechnung Verankerter Hangebr$\ddot{u}$cken,  Berlin, 1935.

%\item{[21']}~ F. Bleich, C. B. McCullough, R. Rosecrans, G. S. Vincent, The mathematical theory of vibration in suspension bridges.  US Government Printing Office, Washington, 1950.

%\item{[18]}~ S. P. Timoshenko, Theory of suspension bridges. I. J. Franklin Inst. 235 (1943), 213-238.

%\item{[19]}~ D. J. Peery, An influence-line analysis for suspension bridge. Translating of ASCE, 1955, 121:463-510.

%\item{[5]}~ L. H. Tsien, A simplified method of analyzing suspension bridges. Trans. ASCE, 114, 1949.

%Irvine, H. M. (1981). Cable structures. MIT University Press, Cambridge, Mass.

\item{[28]}~ B. Semper, Finite element methods for suspension bridge models, Computers. Math. Applic. 25(5) (1993),  77-91.

\item{[29]}~ B. Semper, Finite element approximation of a fourth order integro-differential equation, Appl. Math. Lett. (1) (1994), 59-62.

\item{[30]}~ R. C. Mittal, R. Jiwari, A spectral method for suspension bridge model, International Journal of Applied Mathematics and Mechanics. 5(5) (2009), 66-75.

\item{[31]}~ Q. W. Ren, Q. Q. Zhuang, Legendre-galerkin spectral approximation of a class of fourth-order integro-differential equations, Mathematica Numerica Sinica. 35(2) (2013), 125-136.

\item{[32]}~ E. Aruchunan, Y. Wu, B. Wiwatanapataphee, et al. A new variant of arithmetic mean iterative method for fourth order integro-differential equations solution, 2015 IEEE Third International Conference on Artificial Intelligence, Modelling and Simulation, IEEE. 2015.

\item{[33]}~ A. Zeeshan, M. Atlas, Optimal solution of integro-differential equation of suspension bridge model using Genetic Algorithm and Nelder-Mead method, Journal of the Association of Arab Universities for Basic and Applied Sciences. 24(1) (2017), 310-314.

\item{[34]}~ A. C. Lazer, P. J. Mckenna, Large-amplitude periodic oscillations in suspension bridge: some new connections with nonlinear analysis, SIAM Rev. 32 (1990), 537-578.

\item{[35]}~ R. H. Plaut, F. M. Davis, Sudden lateral asymmetry and torsional oscillations of section models of suspension bridges, J. Sound and Vibration. 307 (2007), 894-905.

\item{[36]}~ F. Gazzola, R. Pavani, Wide oscillations finite time blow up for solutions to nonlinear fourth order differential equations, Arch. Ration. Mech. Anal. 207 (2013), 717-752.

\item{[37]}~ W. Lacarbonara, Nonlinear Structural Mechanics: Theory, Dynamical Phenomena and Modeling[M]. Springer Publishing Company, Incorporated, 2013.  Springer, New York, 2013.

\item{[38]}~ F. Gazzola, M. Jleli, B. Samet, On the melan equation for suspension bridges, J Fixed Point Theory Appl. 16 (2014), 159188.

\item{[39]}~ F. Gazzola, Mathematical models for suspension bridges:
Nonlinear Structural Instability. Springer, New York, 2015.

\item{[40]}~ J. X. Wang, Existence, uniqueness, and approximability of solutions to the classical Melan equation in suspension bridges, Journal of Computational and Applied Mathematics. 485 (2026), 117520.

\item{[41]}~ R. C. Hibbeler, Structural Analysis, 8th Edition. Prentice Hall, 2012.

\item{[42]}~ N. J. Gimsing, C. T. Georgakis, Cable Supported Bridges: Concept and Design, 3rd ed.  John Wiley \& Sons, 2011.

\item{[43]}~ P. B. Baily, L. F. Shampine, P. E. Waltman, Nonlinear two point boundary value problems. New York: Academic Press, 1968.

\item{[44]}~ Q. A. Dang, Q. L. Dang, T. K. Q. Ngo, A novel efficient method for nonlinear boundary value problems. Numer. Algor. 76 (2017), 427-439.

\item{[45]}~ Q. A. Dang, T. K. Q. Ngo, Existence results and iterative method for solving the cantilever beam equation with fully nonlinear
term, Nonlinear Anal. RWA. 36 (2017) 56-68.

\item{[46]}~  Q. A. Dang, T. H. Nguyen, The unique solvability and approximation of BVP for a nonlinear fourth order kirchhoff type equation, East asian. J. Appl. Math. 8(2018), 323-335.

\item{[47]}~ Y. F. Wei, Q. L. Song, Z. B. Bai, Existence and iterative method for some fourth order nonlinear boundary value problems. Appl. Math. Lett. 87 (2019), 101-107.

\item{[48]}~ American Association of State Highway and Transportation Officials, AASHTO LRFD bridge design specifications. 10th Ed., AASHTO, Washington, DC, 2024.
\end{description}
\end{document}